\documentclass{amsart}
\usepackage{amsfonts,amssymb}
\usepackage[labelfont=normal,labelformat=simple]{subcaption}

\usepackage[square,numbers]{natbib}
\usepackage{hyperref}
\usepackage{mathrsfs,bbm}
\usepackage{tikz,tikzscale,tikz-cd}
\usepackage{diagbox}
\usepackage{float}

\newtheorem{theorem}{Theorem}[section]

\theoremstyle{definition}

\newtheorem{definition}[theorem]{Definition}
\newtheorem{example}[theorem]{Example}

\newtheorem{remark}[theorem]{Remark}

\title{Characterization and Topological Behavior of Homomorphism Tree-Shifts}
\author[Jung-Chao Ban]{Jung-Chao Ban}
\address[Jung-Chao Ban]{Department of Mathematical Sciences, National Chengchi University, Taipei 11605, Taiwan, ROC.}
\address{Math. Division, National Center for Theoretical Science, National Taiwan University, Taipei 10617, Taiwan. ROC.}
\email{jcban@nccu.edu.tw}

\author[Chih-Hung Chang]{Chih-Hung Chang}
\address[Chih-Hung Chang]{Department of Applied Mathematics, National University of Kaohsiung, Kaohsiung 81148, Taiwan, ROC.}
\email{chchang@nuk.edu.tw}

\author[Wen-Guei Hu]{Wen-Guei Hu}
\address[Wen-Guei Hu]{College of Mathematics, Sichuan University, Chengdu, 610064, China}
\email{wghu@scu.edu.cn}
\author[Guan-Yu Lai]{Guan-Yu Lai}
\author[Yu-Liang Wu]{Yu-Liang Wu}
\address[Guan-Yu Lai and Yu-Liang Wu]{Department of Applied Mathematics, National Chiao Tung University, Hsinchu 30010, Taiwan, ROC.}
\email{guanyu.am04g@g2.nctu.edu.tw; s92077.am08g@nctu.edu.tw}
\thanks{Ban and Chang are partially supported by the Ministry of Science and Technology, ROC (Contract No MOST 109-2115-M-004-002-MY2 and 109-2115-M-390-003-MY3). Hu is partially supported by the National Natural Science Foundation of China (Grant No.11601355).}

\newcommand\norm[1]{\lvert #1 \rvert}
\newtheorem{namedthm}{Theorem}[section]

\begin{document}
    \maketitle
	\begin{abstract}
	The purpose of this article is twofold. On one hand, we reveal the equivalence of shift of finite type between a one-sided shift $X$ and its associated hom tree-shift $\mathcal{T}_{X}$, as well as the equivalence in the sofic shift. On the other hand, we investigate the interrelationship among the comparable mixing properties on tree-shifts as those on multidimensional shift spaces. They include irreducibility, topologically mixing, block gluing, and strong irreducibility, all of which are defined in the spirit of classical multidimensional shift, complete prefix code (CPC), and uniform CPC. In summary, the mixing properties defined in all three manners coincide for $\mathcal{T}_{X}$. Furthermore, an equivalence between irreducibility on $\mathcal{T}_{A}$ and irreducibility on $X_A$ are seen, and so is one between topologically mixing on $\mathcal{T}_{A}$ and mixing property on $X_A$, where $X_A$ is the one-sided shift space induced by the matrix $A$ and $T_A$ is the associated tree-shift. These equivalences are consistent with the mixing properties on $X$ or $X_A$ when viewed as a degenerate tree-shift.
	\end{abstract}
	
	\section{Introduction}
	
	Let $\mathcal{A}$ be a finite alphabet. A pattern on $\mathbb{Z}^d$ over $\mathcal{A}$ is a function from a finite subset of $\mathbb{Z}^d$ to $\mathcal{A}$. Given a set of patterns $\mathcal{F}$; a shift space $\mathsf{X}_{\mathcal{F}} \subseteq \mathcal{A}^{\mathbb{Z}^d}$ is the set of configurations in which patterns from $\mathcal{F}$ do not appear. Translation of configurations is a natural $\mathbb{Z}^d$ action on $\mathsf{X}_{\mathcal{F}}$ and makes it a dynamical system; $\mathsf{X}_{\mathcal{F}}$ is a shift of finite type (SFT) if $\mathcal{F}$ is finite. The study of SFTs plays an essential role and is rife with numerous undecidability issues whenever $d \geq 2$. For the case where $d=1$, the algorithm deciding the emptiness and the existence of periodic points of an SFT $\mathsf{X}_{\mathcal{F}}$ come immediately from its essential graph representation; $\mathsf{X}_{\mathcal{F}}$ contains dense periodic points if it is irreducible \cite{Kit-1998, LM-1995}. For $d \geq 2$, however, it is not even decidable if $\mathsf{X}_{\mathcal{F}}$ is nonempty; there is an aperiodic SFT with positive topological entropy (cf.~\cite{Berger-MAMS1966, Kari-DM1996, Robinson-IM1971} for instance).
	
	Several topological conditions were introduced to assure certain types of dynamical behavior, such as the density of periodic points, positive topological entropy, and chaos. Boyle \emph{et al.}~\cite{BPS-TAMS2010} demonstrated that every two-dimensional block gluing SFT has dense periodic points while the density of periodic points remains to be open for general $\mathbb{Z}^d$ SFTs, yet for $d \ge 3$, the existence of periodic points remains unknown even for block gluing shifts. Furthermore, every nontrivial block gluing shift space is of positive topological entropy. A possible reason for these differences between one- and multidimensional shift spaces is the spatial structure:
	$\mathbb{Z}$ is a finitely generated group with no relations while $\mathbb{Z}^d$ is not, and this difference also exists in their semigroup counterparts, $\mathbb{Z}_+$ shifts and $\mathbb{Z}_+^d$ shifts.
	
	Shifts on trees (also known as tree-shifts) have received extensive attention in recent years \cite{AB-TCS2012, AB-ToCS2013}. Such shifts exhibit the natural structure of one-dimensional symbolic dynamics while equipped with multiple directional shift transformations. This interesting combination of properties makes the subshifts on trees an intermediate class of symbolic dynamics between one-dimensional shift spaces and multidimensional shift spaces.
	
	Given the richness in interesting properties of multidimensional shift spaces as well as inspired by the physical models, the elucidation of \emph{hom-shifts} is imperative. A hom-shift is a nearest neighbor SFT on $\mathbb{Z}^{d}$ or $\mathbb{Z}_+^{d}$ such that if $a$, $b\in \mathcal{A}$ are forbidden to sit next to each other in some direction, so are they in all coordinate directions. Many important SFTs arise as hom-shifts, for instance, hard square shift. Chandgotia and Marcus \cite{chandgotia2018mixing} studied the mixing properties of hom-shifts and related them to some questions in graph theory therein. Aside from generalization in $\mathbb{Z}^d$ or $\mathbb{Z}_+^d$, the \emph{hom tree-SFT} is also considered as an alternative path on the study of abstract tree-shifts of finite type. The description of hom tree-SFT is as follows. Let $T$ be the free monoid generated by $\Sigma = \{0, 1, \cdots, k-1\}$, $\mathcal{M}_{n}(\{0,1\})$ be the set of $n \times n$ binary matrices, and $A \in \mathcal{M}_{n}(\{0,1\})$ be given. The set
	\[
	\mathcal{T}_{A}:=\{t \in \mathcal{A}^{T}: A(t_{w}, t_{w s})=1 \text{ for } w \in T \text{ and } s \in \Sigma \}
	\]
	is called a \emph{hom tree-SFT}. Mairesse and Marcovici \cite{mairesse2017uniform} considered hom tree-SFTs for $k=2$; they constructed a stationary Markov measure $\mu$ out of adjacency matrix $P$ induced from $A$ and showed that $\mu $ is a Markov uniform measure on the tree. The class of hom tree-SFTs is a particular case of \emph{Markov tree-shifts $\mathcal{T}_{A_0, A_1, \ldots, A_{k-1}}$} considered by Ban and Chang \cite{Ban2017a,Ban2017} when $A_i =A$ for $0 \leq i \leq k-1$, where $\mathcal{T}_{A_0,A_1, \ldots, A_{k-1}}$ is defined as
	\begin{equation} \label{eq:markov_tree_shift}
	    \mathcal{T}_{A_0,A_1, \ldots, A_{k-1}}:=\{t \in \mathcal{A}^{T} : A_i(t_w,t_{wi})=1 \text{ for } w\in T \text{ and } 0 \leq i \leq k-1 \}.
	\end{equation}
	Ban and Chang showed that every tree-shift of finite type is conjugate to a Markov tree-shift; a survey of topological properties is also done therein and is related to chaotic behavior of $\mathcal{T}_{A_0,A_1, \ldots, A_{k-1}}$ \cite{Ban2017a,Ban2017}. 
	
	Similar to illustration of hom tree-SFT, for each one-dimensional shift $X \subset \mathcal{A}^{\mathbb{Z}_+}$ there is an associated \emph{hom tree-shift} $\mathcal{T}_{X}$ defined as
	$$
	\mathcal{T}_{X} := \{t \in \mathcal{A}^T: (t_{w_1 w_2 \ldots w_i})_{i \in \mathbb{Z}_+} \in X \text{ for any sequence } (w_i)_{i \in \mathbb{Z}_+} \text{ in } \Sigma\}.
	$$
	In other words, $\mathcal{T}_{X}$ is the set consists of configurations whose projection on any infinite path is in $X$. Petersen and Salama \cite{PS-TCS2018,Petersen2020entropy} first proposed the class of tree-shifts and demonstrated that the topological entropy of $X$ is no larger than the topological entropy of $\mathcal{T}_{X}$. An immediate result is $X$ being topologically mixing implies $\mathcal{T}_{X}$ is of positive topological entropy. It is noteworthy that this definition is related to the \emph{$d$-fold axial power of $X$ in \cite{Louidor2013} or \emph{isotropic shift} \cite{10.1112/plms/pdu029}, from which the definition of tree-shift can be adapted to
	\begin{equation*}
	    \mathcal{T}_X' := \{t \in \mathcal{A}^T: (t_{g i^n})_{n \in \mathbb{Z}_+} \in X \text{ for any } g \in T, i \in \Sigma\}.
	\end{equation*}
	It is seen from above the definitions of $\mathcal{T}_X$ and $\mathcal{T}_X'$ coincide if $X$ is a \emph{Markov shift}, i.e. $X = X_A$ is defined by a $d \times d$ adjacency matrix $A \in \mathcal{M}_{d}(\{0,1\}) = \{0,1\}^{d \times d}$. Of independent interest, the impact of the nuance in the definitions is demonstrated in Remark \ref{rmk:axial_product}.
	}
	
	This paper aims to investigate the relations of topological properties between $X$ and $\mathcal{T}_{X}$. The upcoming section characterizes two criteria for determining whether $\mathcal{T}_{X}$ is a tree-SFT or a sofic tree-shift as follows.
	\renewcommand\thenamedthm{\ref{thm:SFT_equivalence}}
	\begin{namedthm}
		Suppose $X$ is a shift space. $X$ is an SFT if and only if $\mathcal{T}_{X}$ is a tree-SFT.
	\end{namedthm}
	\renewcommand\thenamedthm{\ref{thm:sofic_equivalence}}
	\begin{namedthm}
		Suppose $Y$ is a shift space. $Y$ is a sofic shift if and only if $\mathcal{T}_{Y}$ is a sofic tree-shift.
	\end{namedthm}
	\noindent It is noteworthy that in the latter theorem above, the covering space of $\mathcal{T}_Y$ is not $\mathcal{T}_X$ in general, for which an example is provided in Example \ref{ex:even_shift}.
	
	Following such a fundamental classification of $\mathcal{T}_{X}$, an extensive elucidation reveals the connections between mixing properties for shift spaces on $\mathbb{Z}^d$ and trees. Ban and Chang considered several mixing properties for tree-shifts through complete prefix codes (CPCs, cf.~\cite{Ban2017a, Ban2017}. Such sets are called complete prefix sets (CPS) in the second paper). The relations of mixing properties, such as strong irreducibility (SI), block gluing (BG), topologically mixing (TM), and irreducibility (IR), in the sense of the classical multidimensional shift, CPC, and uniformly CPC, are delivered in Section 3 as follows.
	\renewcommand\thenamedthm{\ref{thm:Z_top}}
	\begin{namedthm}
		Suppose $\mathcal{T}_{}$ is a tree-shift. Then, the following implications hold:
		\item[1.] $\mathcal{T}_{}$ is SI if $\mathcal{T}_{}$ is CPC USI. 
		\item[2.] $\mathcal{T}_{}$ is BG if $\mathcal{T}_{}$ is CPC UBG. 
		\item[3.] $\mathcal{T}_{}$ is TM if $\mathcal{T}_{}$ is CPC BG. 
		\item[4.] $\mathcal{T}_{}$ is TM if $\mathcal{T}_{}$ is BG. 
	\end{namedthm}
	\noindent Section 3 is further devoted to studying mixing properties for abstract tree-shift $\mathcal{T}_{}$, hom tree-shift $\mathcal{T}_{X}$, and hom tree-SFT $\mathcal{T}_{A}$. Figure \ref{fig:mixing_relation_Z(M)} illustrates the collapse of mixing properties for tree-SFTs. More specifically,
	\begin{figure*}[h]
		\centering
		\includegraphics[]{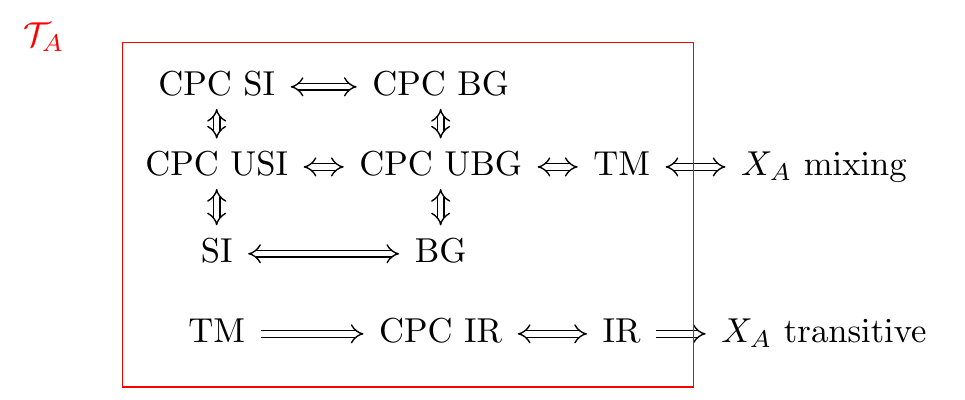}
		\renewcommand{\thefigure}{\ref{fig:mixing_relation_Z(M)}}
		\caption{Relations between mixing properties on Markov tree-shifts $\mathcal{T}_{A}$ and Markov shift spaces $X_A$}
		\setcounter{figure}{0}
	\end{figure*}
	\renewcommand\thenamedthm{\ref{thm:Z(X)_top_with_X}}
	\begin{namedthm}
		Suppose $X$ is a shift space. The following implications hold: 
		\item[1.] $X$ is mixing if $\mathcal{T}_{X}$ is TM.
		\item[2.] $X$ is transitive if $\mathcal{T}_{X}$ is IR.
	\end{namedthm}
	\renewcommand\thenamedthm{\ref{thm:Z(X_M)_top}}
	\begin{namedthm}
		Suppose $A$ is an adjacency matrix. The following implications hold:
		\item[1.] $X_A$ is mixing if and only if $\mathcal{T}_{A}$ is UBG.
		\item[2.] $X_A$ is transitive if and only if $\mathcal{T}_{A}$ is IR.
	\end{namedthm} 
	\section{Characterizations for $\mathcal{T}_{X}$}
	
	This section is devoted to classifying two essential types of tree-shifts. After introducing definitions and notations of tree-shifts, Theorems 2.1 and 2.5 reveal necessary and sufficient conditions for characterizing whether a tree-shift is a shift of finite type or a sofic shift.
	
	\subsection{Notations and Definitions}
	
	Let $k\geq 2$ and $\Sigma_k=\{ 0,1,...,k-1 \}$. A \emph{$k$-tree} $\Sigma_k^\ast = \cup_{n \geq 0} \Sigma_k^n$ is the set of all finite words generated by $\Sigma_k$, where $\Sigma_k^n$ consists of $n$-words and $\Sigma_k^0 = \{\epsilon\}$. The empty word $\epsilon$ denotes the root of the tree and is the only word of zero length. A natural visualization of $\Sigma_k^\ast$ is the Cayley graph of a free semigroup on $k$ generators.
	
	Let $\mathcal{A}=\{0,1,...,d-1\}$ be a finite alphabet. A \emph{labeled tree} is a function $t : \Sigma_k^\ast \rightarrow \mathcal{A}$, and $t_w := t(w)$ denotes the label on the node  $w\in \Sigma_k^\ast$. Let $\mathcal{A}^{\Sigma_k^*}$ be the set of all labeled trees. Define the shift action $\sigma: \Sigma_k^* \times \mathcal{A}^{\Sigma_k^*} \to \mathcal{A}^{\Sigma_k^*}$ as 
	$$
	(\sigma_s t)_w := \sigma (s, t)_w = t_{sw} \quad \text{for all} \quad s, w \in \Sigma_k^*.
	$$
	Set $\Delta_n := \cup_{i=0}^{n}\Sigma_k^i$ to be the initial subtree of the $k$-tree. Notably, $\Delta_n$ is of height $n+1$. An \emph{$n$-block} is a function $u : \Delta_n \rightarrow \mathcal{A}$. In particular, we write every $1$-block $u$ as a $(k+1)$-tuple $(u_\epsilon; u_0, u_1, \cdots, u_{k-1})$. A finite subset $S \subseteq \Sigma_k^\ast$ is said to be \emph{prefix-closed} if each prefix of $S$ lies in $S$. A \emph{pattern} is a function $u: S \to \mathcal{A}$ defined on a finite prefix-closed subset $S$, where $S$ is the \emph{support} (or \emph{shape}) of $u$ and is written as $s(u)=S$. We say that a pattern $u$ appears in a labeled tree $t$ if there is a node $s \in \Sigma_k^\ast$ such that $t_{s w}=u_w$ for all $w \in s(u)$; otherwise, $t$ avoids $u$. A \emph{tree-shift} is a set $\mathcal{T}_{} \subseteq \mathcal{A}^{\Sigma_k^*}$ of labeled trees which avoid all of a certain set of forbidden blocks. A pattern $u$ is admissible in $\mathcal{T}_{}$ if there exists a labeled tree $t \in \mathcal{T}_{}$ such that $(\sigma_s t)|_{s(u)}=u$ for some $s \in \Sigma_k^\ast$; otherwise, $u$ is forbidden in $\mathcal{T}_{}$. Denote by $\mathcal{B}_m(\mathcal{T}_{})$ the set of all admissible $m$-blocks in $\mathcal{T}_{}$; $\mathcal{B}(\mathcal{T}_{}) := \cup_{m \ge 0} \mathcal{B}_m(\mathcal{T}_{})$ refers to the set of all admissible blocks in $\mathcal{T}_{}$.
	
	A subset $\mathbf{s}=\{s_i\}_{i \ge 0} \subset \Sigma_k^\ast$ is called a \emph{chain} if $s_0 = \epsilon$ and $s_{i+1} \in s_i \Sigma_k$ for all $i \ge 0$; in other words, $\mathbf{s}$ is an infinite path initiated at the root. Define a \emph{projection} $\pi_{\mathbf{s}}: \mathcal{A}^{\Sigma_k^*} \to \mathcal{A}^{\mathbb{Z}_+}$ as $(\pi_{\mathbf{s}} t)_i = t_{s_i}$ for $i \in \mathbb{Z}_+$, where $\mathbb{Z}_+:=\mathbb{N} \cup \{0\}$. Suppose $X \subseteq \mathcal{A}^{\mathbb{Z}_+}$ is a shift space. Define the tree-shift $\mathcal{T}_{X}$ by 
	\[
	\mathcal{T}_{X}:=\{t \in \mathcal{A}^{\Sigma_k^\ast}: \pi_{\mathbf{s}} (t) \in X \text{ for every chain } \mathbf{s}\}.
	\]
	In particular, suppose $X \subseteq \mathcal{A}^{\mathbb{Z}_+}$ is a Markov shift by an adjacency matrix $A \in \mathcal{M}_{d}(\{0,1\})=\{0,1\}^{d \times d}$. We denote the tree-shift $\mathcal{T}_{X}$ by $\mathcal{T}_{A}$. Note that every Markov shift is also a \emph{vertex shift}, i.e., $X_A$ has a directed graph representation $\mathsf{G}=(\mathsf{V},\mathsf{E})$ with the vertex set $\mathsf{V}=\mathcal{A}$ and the edge set $\mathsf{E}=\{(\mathsf{v}_i, \mathsf{v}_j) \in \mathsf{V} \times \mathsf{V}: A(\mathsf{v}_i,\mathsf{v}_j)=1\}$. For each $\mathsf{e}= (\mathsf{v}_i, \mathsf{v}_j) \in \mathsf{E}$, we denote the initial and terminal states of $\mathsf{e}$ as $\mathsf{i}(\mathsf{e})=\mathsf{v}_i$ and $\mathsf{t}(\mathsf{e})=\mathsf{v}_j$, respectively. A Markov shift can also be determined by a directed graph, which is inherently equipped with an adjacency matrix (cf.~\cite{LM-1995}), and this definition is an equivalent definition as above. Hence, we also write $X_A=X_{\mathsf{G}}$. For each shift space $X_{\mathsf{G}}$, the \emph{edge shift} of $X_{\mathsf{G}}$ is a Markov shift $X_{\mathsf{G}'}$ with its directed graph representation $\mathsf{G}'=(\mathsf{V}',\mathsf{E}')$, where $\mathsf{V}'=\mathsf{E}$ is its vertex set and $\mathsf{E}'= \{(\mathsf{e}_1,\mathsf{e}_2) \in \mathsf{E} \times \mathsf{E}: \mathsf{t}(\mathsf{e}_1)=\mathsf{i}(\mathsf{e}_2)\}$ is its edge set.
	
	For the sake of convenience, we express $w \in \Sigma_k^n$ as $w=w_0 w_1 \cdots w_n$ in the rest of this elucidation, where $w_0=\epsilon$ and $w_i \in \Sigma_k$. In addition, this paper considers the case $k=2$ for the clarity of discussion, while all of the results generalize to arbitrary $k \in \mathbb{N}$. In the following, we refer to $\Sigma_2$ as $\Sigma$. 
	
	\subsection{Shifts of Finite Type}
	
	A shift space $X = X_{\mathcal{F}}$ is an SFT if the forbidden set $\mathcal{F}$ is finite; $X$ is $m$-step if $\mathcal{F} \subset \mathcal{A}^{m+1}$. A \emph{tree-shift of finite type} is defined in a similar aspect. $\mathcal{T}_{}$ is a tree-SFT if every labeled tree $t \in \mathcal{T}_{}$ avoids any blocks of a finite set, consisting of $n$-blocks for some $n\in \mathbb{N}$. The following theorem reveals a natural and intrinsic characterization of one-dimensional SFTs and tree-SFTs.
	
	\begin{theorem} \label{thm:SFT_equivalence}
		Suppose $X$ is a shift space. $X$ is an SFT if and only if $\mathcal{T}_{X}$ is a tree-SFT.
	\end{theorem}
	
	\begin{proof}
		Suppose $X$ is an SFT with the alphabet $\mathcal{A}$ and forbidden set $\mathcal{F}_X$. Suppose $\mathcal{F}_X \subseteq \mathcal{A}^m$ for some $m \in \mathbb{N}$. We consider a tree-SFT $\mathcal{T}_{}$ defined by the forbidden set 
		\[
		\mathcal{F}_1:=\{u \in \mathcal{A}^{\Delta_{m-1}}: \exists s_1, \cdots, s_{m-1} \in \Sigma, u_{\epsilon} u_{\epsilon s_1} \cdots u_{\epsilon s_1 \cdots s_{m-1}} \in \mathcal{F}_X\},
		\]
		and show that $\mathcal{T}_{}=\mathcal{T}_{X}$. It is clear that $\mathcal{T}_{X} \subseteq \mathcal{T}_{}$ by definition. Suppose $t \in \mathcal{T}_{}$ and $\mathbf{s}=\{s_i\}_{i \ge 0}$ is a chain. By definition of $\mathcal{F}_1$, 
		\[\left(\pi_{\mathbf{s}}(t)\right)_n \left(\pi_{\mathbf{s}}(t)\right)_{n+1} \cdots \left(\pi_{\mathbf{s}}(t)\right)_{n+m-1} = t_{s_n} t_{s_{n+1}} t_{s_{n+m-1}} \notin \mathcal{F}_X \]
		for every $n \ge 0$. Therefore, $\pi_{\mathbf{s}}(x) \in X$ and the desired result follows. 
		
		Conversely, suppose $\mathcal{T}_{X}$ is a tree-SFT with the alphabet $\mathcal{A}$ and forbidden set $\mathcal{F}_{\mathcal{T}_{X}} \subseteq \mathcal{A}^{\Delta_{m-1}}$ for some $m \in \mathbb{N}$. Define an SFT $Y$ by the forbidden set 
		\[\mathcal{F}_2:=\{v \in \mathcal{A}^m: \exists u \in \mathcal{F}_{\mathcal{T}_{X}}, u_w=v_{\norm{w}}, \forall w \in \Delta_{m-1}\}.\] 
		We then show that $Y=X$. Suppose $x=(x_0, x_1, \ldots) \in X$. Define a labeled tree $t \in \mathcal{A}^{\Sigma^\ast}$ as $t_w = x_{\norm{w}}$ for every $w \in \Sigma^\ast$. Then, it follows from definition of $\mathcal{T}_{X}$ that $t \in \mathcal{T}_{X}$. Note that $(\sigma_w t)_{\Delta_{m-1}} \notin \mathcal{F}_{\mathcal{T}_{X}}$ for every $w \in \Sigma^\ast$ implies $x \in Y$ by the definition of $\mathcal{F}_2$. This finishes the proof of $X \subseteq Y$. On the other hand, if $y = (y_0, y_1, \cdots ) \in Y$, we show that $y \in X$. Indeed, if we define a labeled tree $t \in \mathcal{A}^{\Sigma^\ast}$ as $t_w = y_{\norm{w}}$ for every $w \in \Sigma^\ast$, it is clear that $t \in \mathcal{T}_{X}$. Thus, $y = \pi_{\mathbf{0}} (t) \in X$ by definition of $\mathcal{T}_{X}$, where $\mathbf{0} = \{0^i\}_{i \ge 0}$ is a chain.
	\end{proof}
	
	\begin{example}
		We give an example to illustrate the construction of $\mathcal{F}_1$ and $\mathcal{F}_2$ in Theorem \ref{thm:SFT_equivalence}.\\
		\textbf{1.} Let $X \subseteq \{0,1\}^{\mathbb{Z}_+}$ be the golden mean shift, i.e. the forbidden set for $X$ is $\mathcal{F}=\{11\}$. Then, $\mathcal{F}_{1} = \{(u_\epsilon; u_0, u_1) = (1; 1, 0), (1; 0, 1), (1; 1, 1)\}$. It is seen to be the same set as the forbidden set $\mathcal{F}_{\mathcal{T}_{X}}$ of $\mathcal{T}_{X}$.\\
		\textbf{2.} At first glance $\mathcal{F}_2$ seems artificial. However, one may observe that if $v \in \mathcal{A}^m \setminus \mathcal{B}_m(X)$, then there is a corresponding $(m-1)$-block $u \in \mathcal{A}^{\Delta_{m-1}} \setminus \mathcal{B}_{m-1}(\mathcal{T}_{X})$ such that $u_w = v_{\norm{w}}$ for every $w \in \Delta_{m-1}$. For instance, suppose $\mathcal{T}_{X}$ is a hom tree-shift associated with the golden mean shift $X$. The forbidden set of $\mathcal{T}_{X}$ is given by $\mathcal{F}_{\mathcal{T}_{X}}=\mathcal{A}^{\Delta_1} \setminus \mathcal{B}_1(\mathcal{T}_{X})=\{(u_\epsilon; u_0, u_1) = (1; 1, 0), (1; 0, 1), (1; 1, 1)\}$. One may see that $(u_\epsilon; u_0, u_1) = (1; 1, 1) \notin \mathcal{B}_1(\mathcal{T}_{X})$ since $11 \in \mathcal{F}_X$. Also, one can deduce from definition $11 \in \mathcal{F}_{2}$ indeed.
	\end{example}
	
	\subsection{Sofic Shifts}
	
	A one-dimensional shift space $Y \subseteq \mathcal{A}(Y)^{\mathbb{Z}_+}$ is \emph{sofic} if and only if it is a factor of a one-dimensional SFT, i.e., there are an $m$-step SFT $X \subseteq \mathcal{A}(X)^{\mathbb{Z}_+}$ and a block map $f:\mathcal{B}_n(X) \to \mathcal{A}(Y)$ such that the image of the induced map $f_\ast:X \to \mathcal{A}(Y)^{\mathbb{Z}_+}$ is $Y$. And $\mathcal{T}_{} \subseteq \mathcal{A}(\mathcal{T}_{})^{\Sigma^\ast}$ is a sofic tree-shift if and only if there exist an $m$-step tree-SFT $\mathcal{T}_{}' \subseteq \mathcal{A}(\mathcal{T}_{}')^{\Sigma^\ast}$ and a block map $g:\mathcal{B}_n(\mathcal{T}_{}') \to \mathcal{A}(\mathcal{T}_{})$ such that the image of the induced map $f_\ast:\mathcal{T}_{}' \to \mathcal{A}(\mathcal{T}_{})^{\Sigma^\ast}$ is $\mathcal{T}_{}$. In the following we discuss the relation between sofic $\mathcal{T}_{X}$ and sofic $X$. We note that $\mathcal{B}_1(X) = \mathcal{A}(X) = \mathcal{A}(\mathcal{T}_{X}) = \mathcal{B}_0 (\mathcal{T}_{X})$ are all equal by definition, and thus the following two terms are coined to unify the names of such maps on $X$ and $\mathcal{T}_{X}$. 
	
	\begin{definition}
		Let $X$ be a shift space, and $\mathcal{A}$ be some alphabet. A $1$-block map on $X$, $f: \mathcal{B}_1(X) \to \mathcal{A}$, is called a \emph{symbol map on $X$}, and a $0$-block map $g:\mathcal{B}_0(\mathcal{T}_{X}) \to \mathcal{A}$ is called a \emph{symbol map on $\mathcal{T}_{X}$}.
	\end{definition}
	
	Note that in the above definition, the 1-sliding block code induced by $f$ is a map $f_\ast: X \to \mathcal{A}_2^{\mathbb{Z}_+}$ defined as $((f_\ast(x))_i)_{i \ge 0} = (f(x_i))_{i \ge 0}$, and the 0-sliding block code induced by $g$ is a map $g_\ast: \mathcal{T}_{X} \to \mathcal{A}_2^{\Sigma_k^\ast}$ defined as $((g_\ast(x))_w)_{w \in \Sigma_k^\ast} = (g(x_w))_{w \in \Sigma_k^\ast}$. Before we state the main result, we need the following useful theorem.
	
	\begin{theorem} \label{thm:sliding_block_code}
		$\mathcal{T}_{}$ is a sofic tree-shift if and only if there exist a Markov tree-shift $\mathcal{T}_{}'$ (see \eqref{eq:markov_tree_shift}) and a symbol map $g:\mathcal{A}(\mathcal{T}_{}') \to \mathcal{A}(\mathcal{T}_{})$ such that the image of $g_\ast:\mathcal{T}_{}' \to \mathcal{A}(\mathcal{T}_{})^{\Sigma^\ast}$ is $\mathcal{T}_{}$.
	\end{theorem}
	
	\begin{proof}
		By definition of sofic tree-shift, there exist an $m$-step tree-SFT $\mathcal{T}_{}''$ and an $n$-block map $g'':\mathcal{B}_n (\mathcal{T}_{}'') \to \mathcal{A}(\mathcal{T}_{})$ such that $g''_\ast (\mathcal{T}_{}'') = \mathcal{T}_{}$. We construct a Markov tree-shift $\mathcal{T}_{}'$ as follows. Put $N = \max\{m,n\}$. It is shown in \cite[Proposition~2.4]{Ban2017} that there exist a Markov tree-shift $\mathcal{T}_{}'=\mathcal{T}_{}''^{[N]}$ conjugate to $\mathcal{T}_{}''$, and a bijective $N$-block map $g':\mathcal{B}_{N}(\mathcal{T}_{}'') \to \mathcal{A}(\mathcal{T}_{}')$ such that the induced sliding block code $g'_\ast$ satisfies $g'_\ast (\mathcal{T}_{}'')=\mathcal{T}_{}'$. By defining the symbol map $g(a)=g''({g'}^{-1}(a)|_{\Delta_n})$ on $\mathcal{T}_{}'$, it is obvious $g_\ast (\mathcal{T}_{}')=\mathcal{T}_{}$.
	\end{proof}
	
	The following theorem provides a natural and intrinsic characterization of one-dimensional sofic shifts and sofic tree-shifts. Notably, the covering space of $\mathcal{T}_Y$ needs not be $\mathcal{T}_X$ generally; the construction of the proper covering space is the main difficulty. See Example 2.6 for instance.
	
	\begin{theorem} \label{thm:sofic_equivalence}
		Suppose $Y$ is a shift space. $Y$ is a sofic shift if and only if $\mathcal{T}_{Y}$ is a sofic tree-shift.
	\end{theorem}
	\begin{proof}
		\textbf{1}. We first prove that $Y$ is a sofic shift if $\mathcal{T}_{Y}$ is a sofic
		tree-shift. Since $\mathcal{T}_{Y}$ is a sofic tree-shift, by Theorem \ref{thm:sliding_block_code} there is a Markov tree-shift $\mathcal{T}_{}$ defined as 
		\[
		\mathcal{T}_{}=\mathcal{T}_{A_{0},A_{1}}:=\{t\in \mathcal{A}(\mathcal{T}_{})^{\Sigma
			^{\ast }}:A_{i}(t_{w},t_{wi})=1, \text{ for }i=0,1\}\text{,}
		\]%
		and a symbol map, say $g:\mathcal{A}(\mathcal{T}_{})%
		\mathcal{\rightarrow A}(\mathcal{T}_{Y})$, such that $g_{\ast }(\mathcal{T}_{})=\mathcal{T}_{Y}$, in which we assume that $A_{0}$ and $A_{1}\in \mathcal{M}_{d}(\{0,1\}),$
		i.e., $\mathcal{A}(\mathcal{T}_{})=\{0,1, \cdots, d-1\}$, for some $d \in \mathbb{N}$. In the
		following we will construct an SFT $X\subseteq \mathcal{A}(X)^{\mathbb{Z}_+}$ and a symbol map $f:\mathcal{A}(X)\rightarrow \mathcal{A}(Y)$
		such that $f_{\ast }(X)=Y$. From this we deduce that $Y$ is a sofic shift. 
		
		Let $\mathcal{A}(X)=\mathcal{A}(\mathcal{T}_{})$. Define the $d\times d$ matrix $A$ as the coordinatewise maximum of $A_0$ and $A_1$, and define a symbol map $f:\mathcal{A}(X)\rightarrow \mathcal{A}(Y)$ by
		letting $f(a)=g(a)$ for any $a\in \mathcal{A}(X)$. Since $\mathcal{A}(X)=\mathcal{A}(\mathcal{T})$, it is evident that $f$ is
		well-defined. Finally, we define $X=X_{A}$, namely, the Markov shift associated with adjacency matrix $A$. Note that 
		\begin{equation}
		X_{A_{i}}\subseteq X_{A}\text{ for }i=0,1.  \label{3}
		\end{equation}
		
		We first claim that for any chain $\mathbf{s}=\{s_{i}\}_{i\geq 0}$ and
		for any $x\in \mathcal{T}_{}$ we have 
		\begin{equation}
		f_{\ast }(\pi_{\mathbf{s}}(x))=\pi_{\mathbf{s}}(g_{\ast }(x)).  \label{1}
		\end{equation}%
		Indeed, one may verify that
		\[
		\pi_{\mathbf{s}}(g_{\ast }(x))=\pi_{\mathbf{s}}(g_{\ast }((x_{w})_{w\in \Sigma
			^{\ast }}))=\pi_{\mathbf{s}}((g(x_{w}))_{w\in \Sigma ^{\ast
		}})=(f(x_{s_{i}}))_{i\geq 0}\text{,}
		\]%
		and that 
		\[
		f_{\ast }(\pi_{\mathbf{s}}(x))=f_{\ast }\left( (x_{s_{i}})_{i\geq 0}\right)
		=(f(x_{s_{i}}))_{i\geq 0}=(f(x_{s_{i}}))_{i\geq 0}\text{,}
		\]%
		which shows that the equality (\ref{1}) holds. 
		
		We next claim that $f_{\ast }(X)=Y$. Let $x=(x_{0},x_{1},\cdots )\in X$. From the definition of $A$, we choose a chain $\mathbf{s}=\{s_{i}\}_{i\geq 0}$ according to $x$ such that $s_{i}=s_{0}' s_{1}' \cdots s_{i}'$ and $A_{s_{i}'}(x_{i},x_{i+1})=1$ for all $i\geq 1$.
		Thus there exists $t \in \mathcal{T}_{}$ such that $\pi_{\mathbf{s}}(t)=x$ and $\pi_{\mathbf{s}}(g_{\ast }(t))\in Y$, and we
		conclude that $f_{\ast}(x)=f_{\ast} \left(\pi_{\mathbf{s}}(\mathcal{T}_{})\right) \in Y$. This proves $f_{\ast }(X)\subseteq Y$ . For the converse,
		let $\mathbf{0}=\{0^{i}\}_{i\geq 0}$ be the leftmost chain of $\Sigma ^{\ast }$. Since $g_\ast (\mathcal{T}_{}) = \mathcal{T}_{Y}$, it can be verified by definition that $Y=\pi_{\mathbf{0}}(g_{\ast }(\mathcal{T}_{}))$. By combining this with \eqref{3} and \eqref{1}, it follows
		\[
		Y=\pi_{\mathbf{0}}(g_{\ast }(\mathcal{T}_{}))=f_\ast(\pi_{\mathbf{0}}(\mathcal{T}_{}))=f_\ast (X_{A_0}) \subseteq f_{\ast }(X)\text{.}
		\]%
		Thus $f_{\ast }(X)=Y$.
		
		\textbf{2}. Suppose $Y$ is a sofic shift and $(\mathsf{G},f)$ is the
		corresponding graph representation, where $\mathsf{G}=(\mathsf{V},\mathsf{E})$ is a directed graph and $%
		f:\mathsf{E}\rightarrow \mathcal{A}(Y)$ is a map from the edge set $\mathsf{E}$ to the symbol
		set $\mathcal{A}(Y)$. Note that the vertex shift $X_{\mathsf{G}}$ is conjugate to its edge
		shift (see \cite[Theorem~2.4.10]{LM-1995}), i.e., $X^{(1)}=X_{\mathsf{G}^{(1)}}$, where $%
		\mathsf{G}^{(1)}=(\mathsf{V}^{(1)},\mathsf{E}^{(1)})$, where $\mathsf{V}^{(1)}=\mathsf{E}$ and 
		\[
		\mathsf{E}^{(1)}=\{(\mathsf{e}_{1}, \mathsf{e}_{2}) \in \mathsf{E}\times \mathsf{E}:\mathsf{t}(\mathsf{e}_{1})=\mathbf{i}(\mathsf{e}_{2})\}\text{.} 
		\]%
		By defining a symbol map $f^{(1)}:\mathcal{A}(X^{(1)})\rightarrow \mathcal{A}%
		(Y)$ by $f^{(1)}=f$, it is clear that $f_{\ast }^{(1)}(X^{(1)})=Y$.
		
		Our goal is to construct a subshift of finite type $X^{(2)}=X_{\mathsf{G}^{(2)}},$
		where $\mathsf{G}^{(2)}=(\mathsf{V}^{(2)},\mathsf{E}^{(2)})$, and symbol maps $f^{(2)}:\mathcal{A}%
		(X^{(2)})\rightarrow \mathcal{A}(Y)$, $g^{(2)}:\mathcal{A}%
		(\mathcal{T}_{X^{(2)}})\rightarrow \mathcal{A}(\mathcal{T}_{Y})$ such that 
		\begin{equation}
		f_{\ast }^{(2)}(X^{(2)})=Y\text{ and }g_{\ast }^{(2)}(\mathcal{T}_{X^{(2)}})=\mathcal{T}_{Y}\text{.%
		}  \label{10}
		\end{equation}
		This shows that $\mathcal{T}_{Y}$ is a sofic tree-shift since $\mathcal{T}_{X^{(2)}}$ is a tree-SFT as shown in Theorem \ref{thm:SFT_equivalence}.
		
		The construction will be established as follows. For $\mathsf{v}^{(1)},\mathsf{v}'^{(1)}\in \mathsf{V}^{(1)}=\mathcal{A}(X^{(1)})$, denote by $%
		\mathsf{v}^{(1)}R\mathsf{v}'^{(1)}$ the relation on $\mathsf{v}^{(1)}$ and $\mathsf{v}'^{(1)}$ such that 
		\begin{equation}
		f^{(1)}(\mathsf{v}^{(1)})=f^{(1)}(\mathsf{v}'^{(1)}).  \label{4}
		\end{equation}
		One may verify that it is an equivalence relation. We thus denote by $\{[\mathsf{v}^{(1)}]:\mathsf{v}^{(1)}\in \mathsf{V}^{(1)}\}$ the associated equivalence class.
		Let us introduce a new vertex set 
		\[
		\mathsf{V}^{(2)}=\bigcup_{\mathsf{v}^{(1)}\in \mathcal{A}(X^{(1)})}2^{[\mathsf{v}^{(1)}]}\backslash
		\{\emptyset \}, 
		\]%
		where $2^{[\mathsf{v}^{(1)}]}$ stands for the power set of the $[\mathsf{v}^{(1)}]$. The
		associated edge set $\mathsf{E}^{(2)}$ is defined as follows. For $\mathsf{v}^{(2)}, \mathsf{v}'^{(2)} \subseteq \mathsf{V}^{(1)}$, we say $(\mathsf{v}^{(2)}, \mathsf{v}'^{(2)}) \in \mathsf{E}^{(2)}$ if for every $\mathsf{v}'^{(1)} \in \mathsf{v}'^{(2)}$ there exists $\mathsf{v}^{(1)} \in \mathsf{v}^{(2)}$ such that $(\mathsf{v}^{(1)},  \mathsf{v}'^{(1)}) \in \mathsf{E}^{(1)}$. Let $X=X_{\mathsf{G}^{(2)}}$ be the
		associated SFT, and $f^{(2)}:\mathcal{A}(X^{(2)})\rightarrow \mathcal{A}(Y)$
		be the symbol map which is defined by $f^{(2)}(a^{(2)})=f^{(1)}(a^{(1)})$
		for all $a^{(1)}\in a^{(2)}$. The definition of $f^{(2)}$ is well-defined
		since if $a^{(1)}$ and $b^{(1)}$ belong to the same equivalence class, then
		they have the same image under $f^{(1)}$ by (\ref{4}). We note that there
		is an associated symbol map, say $g^{(2)}:\mathcal{A}(\mathcal{T}_{X^{(2)}})%
		\rightarrow \mathcal{A}(\mathcal{T}_{Y})$ with $g^{(2)}(a^{(2)})=f^{(2)}(a^{(2)})$
		meanwhile. The following two properties are essential for the construction
		of $X^{(2)}$ and $\mathcal{T}_{X^{(2)}}$.
		
		\textbf{i}. We first show that if $\mathsf{v}_{1}^{(2)}\mathsf{v}_{2}^{(2)}\cdots \mathsf{v}_{n}^{(2)}$
		is a path in $\mathsf{G}^{(2)}$, then for all $\mathsf{v}_n^{(1)} \in \mathsf{v}_n^{(2)}$ there exists a path, say $\mathsf{v}_{1}^{(1)}\mathsf{v}_{2}^{(1)}%
		\cdots \mathsf{v}_{n}^{(1)}$ in $\mathsf{G}^{(1)}$ such that $\mathsf{v}_{i}^{(1)}\in $ $\mathsf{v}_{i}^{(2)}$
		for $1\leq i\leq n$. We prove it by induction on $n$. The case when $n=1$
		holds by the definition. Assume the claim holds for $n=N$. Suppose $%
		\mathsf{v}_{1}^{(2)}\mathsf{v}_{2}^{(2)}\cdots \mathsf{v}_{N}^{(2)}\mathsf{v}_{N+1}^{(2)}$ is a path in $\mathsf{G}^{(2)}$%
		, then for each $\mathsf{v}_{N+1}^{(1)}\in \mathsf{v}_{N+1}^{(2)}$ there exists, by definition of $\mathsf{E}^{(2)}$, an edge $(\mathsf{v}_{N}^{(1)}, \mathsf{v}_{N+1}^{(1)})$ such that $\mathsf{v}_{N}^{(1)}\in \mathsf{v}_{N}^{(2)}$. By
		induction, there exists a path $\mathsf{v}_{1}^{(1)}\mathsf{v}_{2}^{(1)}\cdots \mathsf{v}_{N}^{(1)}$ in 
		$\mathsf{G}^{(1)}$ such that $\mathsf{v}_{i}^{(1)}\in \mathsf{v}_{i}^{(2)}$ for $i\leq N$. Thus $%
		\mathsf{v}_{1}^{(1)}\mathsf{v}_{2}^{(1)}\cdots \mathsf{v}_{N}^{(1)}\mathsf{v}_{N+1}^{(1)}$ is also a path in $%
		\mathsf{G}^{(1)}$ such that $\mathsf{v}_{i}^{(1)}\in \mathsf{v}_{i}^{(2)}$ for $i\leq N+1$.
		
		\textbf{ii}. If $\mathsf{v}_{1}^{(2)} \mathsf{v}_{1}^{(2)} \cdots \mathsf{v}_{n}^{(2)}$ and $\mathsf{v}_{1}'^{(2)} \mathsf{v}_{1}'^{(2)} \cdots \mathsf{v}_{n}'^{(2)}$ are paths in $\mathsf{G}^{(2)}$ with $f^{(2)}(\mathsf{v}_{i}^{(2)}) = f^{(2)}(\mathsf{v}_{i}'^{(2)})$ for every $1 \le i \le n$, then by writing $\mathsf{v}_{1}''^{(2)} = \mathsf{v}_{1}^{(2)} \cup \mathsf{v}_{1}'^{(2)}$, we get $\mathsf{v}_{1}''^{(2)} \mathsf{v}_{1}''^{(2)} \cdots \mathsf{v}_{n}''^{(2)}$ is a path in $\mathsf{G}^{(2)}$. Indeed, it can be verified by definition of $\mathsf{E}^{(2)}$ that $(\mathsf{v}_{i}''^{(2)}, \mathsf{v}_{i+1}''^{(2)} ) \in \mathsf{E}^{(2)}$ for every $1 \le i \le n-1$.
		
		Now we are ready to show \eqref{10}. It follows from \textbf{i} that
		\begin{equation}
		f_{\ast }^{(2)}(X^{(2)}) \subseteq f_{\ast }^{(1)}(X^{(1)}) = Y.
		\label{7-1}
		\end{equation}%
		Furthermore, one may observe that every path $\mathsf{v}_{1}^{(1)} \mathsf{v}_{2}^{(1)} \cdots \mathsf{v}_{n}^{(1)}$ in $\mathsf{G}^{(1)}$ is also a path in $\mathsf{G}^{(2)}$ by definition of $\mathsf{E}^{(2)}$. Hence, 
		\begin{equation}
		Y=f_{\ast }^{(1)}(X^{(1)})\subseteq f_{\ast }^{(2)}(X^{(2)}).
		\label{7}
		\end{equation}%
		Applying \eqref{7-1}, for every labeled tree $t \in \mathcal{T}_{X^{(2)}}$ and every chain $\mathbf{s}$, we assert that $\pi_{\mathbf{s}}(g_\ast^{(2)}(t)) = f_\ast^{(2)}(\pi_{\mathbf{s}}(t)) \in Y$, since $\pi_{\mathbf{s}}(t) \in X^{(2)}$ by definition. Thus,
		\begin{equation}
		g_{\ast }^{(2)}(\mathcal{T}_{X^{(2)}})\subseteq \mathcal{T}_{Y}\text{.}  \label{8}
		\end{equation}%
		Our next goal is to claim that if $t \in
		\mathcal{T}_{Y}$, there exists $\overline{t} \in \mathcal{T}_{X^{(2)}}$ such that $g^{(2)} (\overline{t}) = t$. It is equivalent to show that for every $n \in \mathbb{N}$, there exists $\overline{u} \in \mathcal{B}_n (\mathcal{T}_{X^{(2)}})$ such that $(g^{(2)} (\overline{u}_{w}))_{w \in \Delta_n} = (t_w)_{w \in \Delta_n}$. Once this claim holds, there exists $\overline{t}^{(n)} \in \mathcal{T}_{X^{(2)}}$ such that $(g^{(2)}(\overline{t}^{(n)}_{w}))_{w \in \Delta_n} = (t_w)_{w \in \Delta_n}$. By compactness of $\mathcal{T}_{X^{(2)}}$, there exists a convergent subsequence $\overline{t}^{(n_k)}$ of $\overline{t}^{(n)}$ with $\lim_{k \to \infty} \overline{t}^{(n_k)} = \overline{t} \in \mathcal{T}_{X^{(2)}}$. Since $g^{(2)}$ is a symbol map, $g^{(2)}_\ast (\overline{t}) = t$, which proves
		\begin{equation}
		\mathcal{T}_{Y} \subseteq g_{\ast }^{(2)}(\mathcal{T}_{X^{(2)}})\text{.} 
		\label{8-1}
		\end{equation}
		By combining \eqref{7-1}, \eqref{7}, \eqref{8} and \eqref{8-1}, we finish the proof of \eqref{10}.
		
		We now show the claim, and we prove it by induction on $n$. The case $n = 1$ holds due to the following reason. If $u=(u_\epsilon; u_0, u_1) \in \mathcal{B}_1(\mathcal{T}_{Y})$, there exist paths $\mathsf{v}^{[0]}_0 \mathsf{v}^{[0]}_1$ and $\mathsf{v}^{[1]}_0 \mathsf{v}^{[1]}_1$ in $\mathsf{G}^{(2)}$ such that 
		\begin{eqnarray}
		f^{(1)}(\mathsf{v}_{0}^{[0]})f^{(1)}(\mathsf{v}_{1}^{[0]})= u_\epsilon u_0, \\
		f^{(1)}(\mathsf{v}_{0}^{[1]})f^{(1)}(\mathsf{v}_{1}^{[1]})= u_\epsilon  u_1.
		\end{eqnarray}%
		By \textbf{ii}, there is an admissible block 
		\[
		\overline{u}=(\overline{u}_\epsilon; \overline{u}_0,\overline{u}_1) = (\mathsf{v}_{0}^{[0]} \cup \mathsf{v}_{0}^{[1]}; \mathsf{v}_{1}^{[0]},\mathsf{v}_{1}^{[1]}) \in \mathcal{B}_1 (\mathcal{T}_{X^{(2)}}),
		\]
		such that 
		\[(g^{(2)}(\overline{u}_\epsilon);g^{(2)}(\overline{u}_0),g^{(2)}(\overline{u}_1)) = (u_\epsilon; u_0, u_1).\]
		Suppose the hypothesis holds for $n=N$. Then, for every $u \in \mathcal{B}_{N+1} (\mathcal{T}_{Y})$, there exists $\overline{u}' \in \mathcal{B}_{N} (\mathcal{T}_{X^{(2)}})$ such that $g^{(2)} (\overline{u}'_w)_{w \in \Delta_N} = (u_w)_{w \in \Delta_N}$. Then, for each $z = \epsilon z_1 z_2 \cdots z_N \in \Sigma^N$, there exist paths $\mathsf{v}_0^{[z 0]} \mathsf{v}_1^{[z 0]} \cdots \mathsf{v}_{n+1}^{[z 0]}$ and $\mathsf{v}_0^{[z 1]} \mathsf{v}_1^{[z 1]} \cdots \mathsf{v}_{n+1}^{[z 1]}$ such that
		\begin{eqnarray}
		f^{(1)}(\mathsf{v}_{0}^{[z 0]})f^{(1)}(\mathsf{v}_{1}^{[z 0]})\cdots f^{(1)}(\mathsf{v}_{n+1}^{[z 0]})
		&=&u_{\epsilon}u_{\epsilon z_{1}}\cdots u_{\epsilon z_{1} \cdots z_{n}0},  \label{5} \\
		f^{(1)}(\mathsf{v}_{0}^{[z 1]})f^{(1)}(\mathsf{v}_{1}^{[z 1]})\cdots f^{(1)}(\mathsf{v}_{n+1}^{[z 1]})
		&=&u_{\epsilon}u_{\epsilon z_{1}}\cdots u_{\epsilon z_{1}\cdots z_{n}1}.  \label{6}
		\end{eqnarray}%
		Hence, applying \textbf{ii}, we extend $\overline{u}'$ to an admissible pattern $\overline{u}''$ in $\mathcal{T}_{X^{(2)}}$ with support $\Delta_N \cup \{z 0, z 1\}$, which is defined as
		\begin{equation*}
		\overline{u}''_w := \begin{cases}
		\overline{u}'_w \cup \mathsf{v}_{\norm{w}}^{[z 0]} \cup \mathsf{v}_{\norm{w}}^{[z 1]}, & \text{ if } w \text{ is a subword of } z;\\
		\mathsf{v}_{\norm{w}}^{[z 0]}, & \text{ if } w = z 0;\\
		\mathsf{v}_{\norm{w}}^{[z 1]}, & \text{ if } w = z 1;\\
		\overline{u}'_w & \text{ otherwise};\\
		\end{cases}
		\end{equation*}
		so that
		\[(g^{(2)}(\overline{u}''_w))_{w \in s(\overline{u}'')} = (u_w)_{w \in s(\overline{u}'')}.\]
		The desired $\overline{u} \in \mathcal{B}_{N+1} (\mathcal{T}_{X^{(2)}})$ can be obtained by repeating the process above to every path $z \in \Sigma^N$ from the leftmost to the rightmost path, and the claim thus holds for all $n \in \mathbb{N}$ by induction.
	\end{proof}
	
	\begin{figure}
		\centering
		\begin{subfigure}[b]{0.45\textwidth}
			\centering
			\includegraphics[]{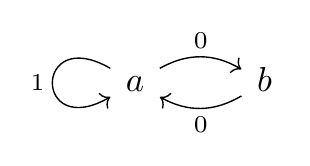}
			\caption{Graph representation $\mathsf{G} = (\mathsf{V},\mathsf{E})$ of cover $X_{\mathsf{G}}$ for even shift $Y$.}
			\label{fig:golden_mean_factor}
		\end{subfigure} \hspace{1em}
		\begin{subfigure}[b]{0.5\textwidth}
			\centering
			\includegraphics[]{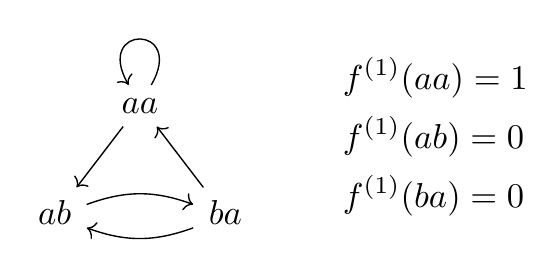}
			\caption{Graph representation $\mathsf{G}^{(1)}=(\mathsf{V}^{(1)},\mathsf{E}^{(1)})$ of the edge shift of $X_{\mathsf{G}}$.}
			\label{fig:golden_mean_edge_shift}
		\end{subfigure}\\
		\begin{subfigure}[b]{0.6\textwidth}
			\centering
			\includegraphics{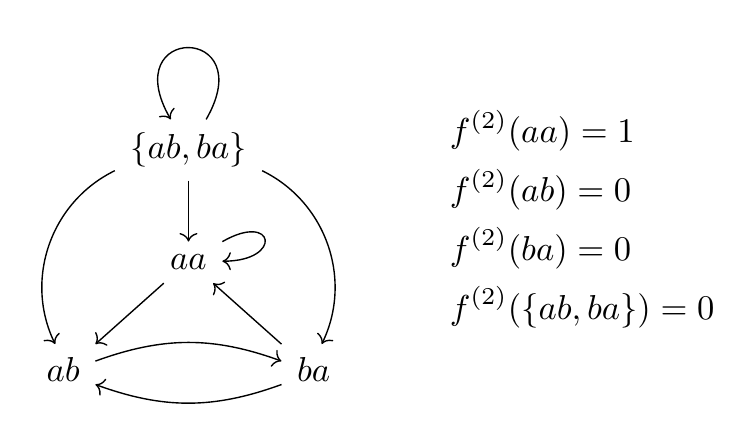}
			\caption{Derived Graph representation $\mathsf{G}^{(2)}=(\mathsf{V}^{(2)},\mathsf{E}^{(2)})$ for the even shift $Y$.}
			\label{fig:even_shift_full_graph}
		\end{subfigure}
		\caption{Construction of the cover for even shift using Theorem \ref{thm:sofic_equivalence}}
	\end{figure}
	
	\begin{figure}
		\resizebox{\textwidth}{!}{%
			\includegraphics[]{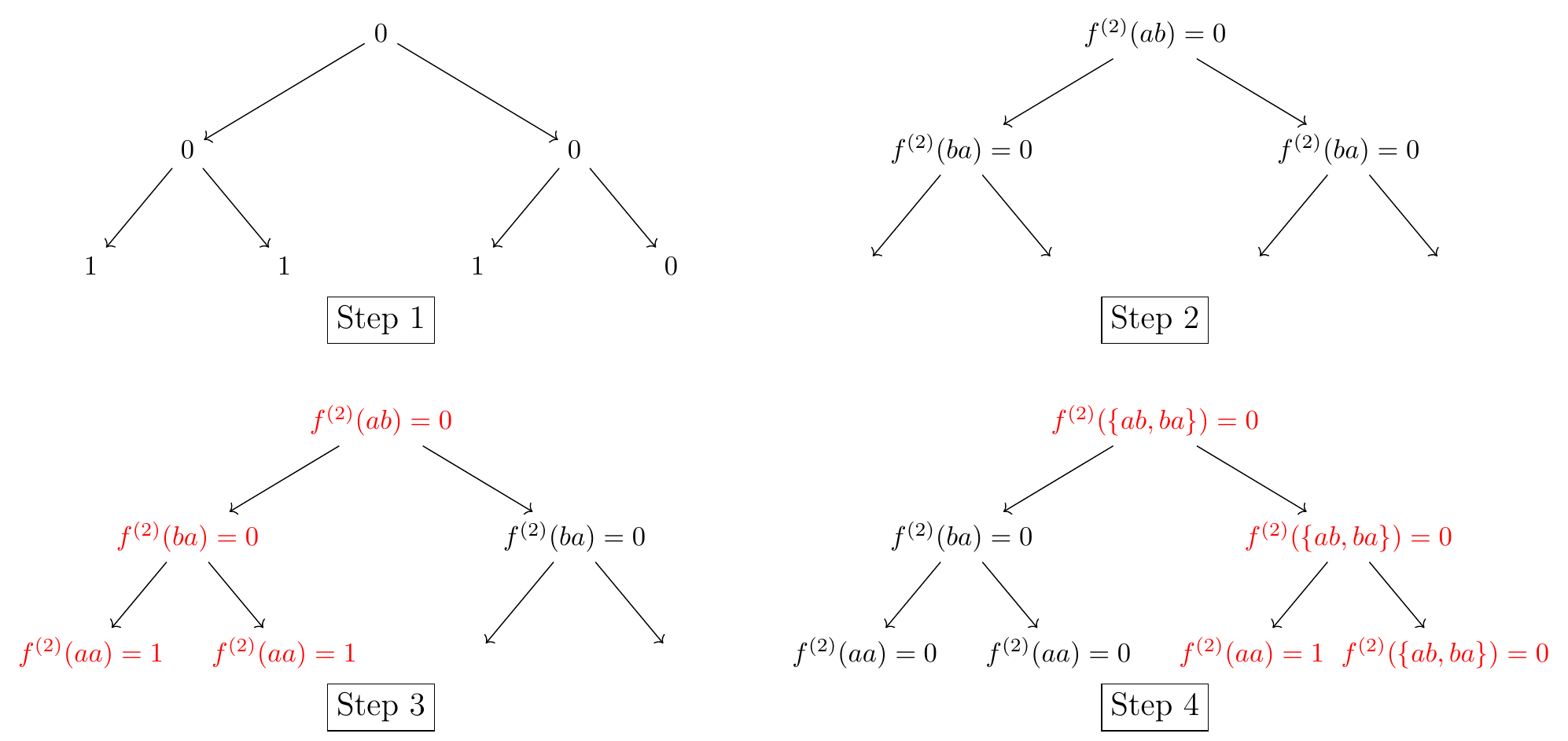}
		}
		\caption{Demonstration of extension process in Theorem \ref{thm:sofic_equivalence}}
		\label{fig:pattern_extension}
	\end{figure}
	\begin{example} \label{ex:even_shift}
		It is known that the even shift $Y$ is a factor of the golden mean shift with graph representation shown in Figure \ref{fig:golden_mean_factor} and the conjugate edge shift $X^{(1)}$ in Figure \ref{fig:golden_mean_edge_shift}. By Theorem \ref{thm:sofic_equivalence}, there exists a cover $X^{(2)}$ by the symbol map $f^{(2)}$ of $Y$ such that $\mathcal{T}_{X^{(2)}}$ is also a cover of $\mathcal{T}_{Y}$ by the symbol map $g^{(2)} = f^{(2)}$, as is shown in Figure \ref{fig:even_shift_full_graph}.  
		
		\textbf{1}. In the proof of necessity of Theorem \ref{thm:sofic_equivalence}, we construct an adjacency matrix $A$ from adjacency matrices $A_0$ and $A_1$ while preserving the image of $f_\ast$. This step is essential since $A_0$ may not coincide with $A_1$ in general. For example, one may consider even shift $Y \subseteq \{0,1\}^{\mathbb{Z}_+}$ and a cover $\mathcal{T}_{A_0,A_1}$ for $\mathcal{T}_{Y}$ as follows. Let $\mathcal{T}_{A_0,A_1}\subseteq\{0,1,2,3,4\}^{\Sigma^\ast}$ be a Markov tree-shift, where
		\begin{equation}
		A_0=
		\left[\begin{array}{c|cc|c|c}
		1 & 1 & 1 & 1 & 1 \\
		\hline
		0 & 1 & 0 & 1 & 0 \\
		0 & 0 & 1 & 1 & 0 \\
		\hline
		0 & 0 & 0 & 0 & 1 \\
		\hline
		0 & 1 & 1 & 1 & 0 \\
		\end{array}\right],
		A_1=
		\left[\begin{array}{c|cc|c|c}
		1 & 1 & 1 & 1 & 1 \\
		\hline
		0 & 1 & 1 & 1 & 0 \\
		0 & 1 & 1 & 1 & 0 \\
		\hline
		0 & 0 & 0 & 0 & 1 \\
		\hline
		0 & 1 & 1 & 1 & 0 \\
		\end{array}\right].
		\end{equation}
		Suppose $g':\{0,1,2,3,4\} \to \mathcal{A}(X^{(2)})$ is a symbol map defined as 
		\begin{equation*}
		\begin{cases}
		g'(0) = \{ab, ba\},\\
		g'(1) = g'(2) = \{aa\},\\
		g'(3) = \{ab\},\\
		g'(4) = \{ba\},
		\end{cases}
		\end{equation*}
		it can be shown that $g'_\ast(\mathcal{T}_{A_0,A_1})=\mathcal{T}_{X^{(2)}}$ by noting that \textbf{a}.~for every $(u'_\epsilon; u'_0, u'_1) \in \mathcal{B}_1(\mathcal{T}_{A_0,A_1})$, $(g'(u'_\epsilon); g'(u'_0), g'(u'_1)) \in \mathcal{B}_1(X^{(2)})$, and \textbf{b}.~by removing symbol 2, it coincides with $X^{(2)}$. As a result, the symbol map $g=g^{(2)} \circ g'$ on $\mathcal{T}_{A_0,A_1}$ satisfies $g_\ast(\mathcal{T}_{A_0,A_1})=\mathcal{T}_{Y}$. On the other hand, by defining
		\begin{equation}
		A:=
		\left[\begin{array}{c|cc|c|c}
		1 & 1 & 1 & 1 & 1 \\
		\hline
		0 & {1} & {1} & 1 & 0 \\
		0 & {1} & {1} & 1 & 0 \\
		\hline
		0 & 0 & 0 & 0 & 1 \\
		\hline
		0 & 1 & 1 & 1 & 0 \\
		\end{array}\right]
		\end{equation}
		as in Theorem \ref{thm:sofic_equivalence}, it is seen that $g_\ast(\mathcal{T}_{A})=\mathcal{T}_{Y}$ by adopting the same argument as above. {It is noteworthy that $A=A_1$ in this case.}
		
		\textbf{2}. In the proof of sufficiency in Theorem \ref{thm:sofic_equivalence}, we construct $\overline{u} \in \mathcal{B}_{N+1}(\mathcal{T}_{X^{(2)}})$ by an extension process. {We now demonstrate the process by our example of even shift, which is illustrated step by step in Figure \ref{fig:pattern_extension} and explained as follows.} Given an admissible $2$-block $u \in \mathcal{B}_2(\mathcal{T}_{Y})$ as in Step 1, we may find an admissible $1$-block $\overline{u}'$ in Step 2. Then, for $00, 01, 10, 11 \in \Sigma^2$, choose the following paths in $\mathsf{G}^{(2)}$.
		\begin{eqnarray*}
		\mathsf{v}_0^{[00]} \mathsf{v}_1^{[00]} \mathsf{v}_{2}^{[00]}=\{ab\}\{ba\}\{aa\}, \\
		\mathsf{v}_0^{[01]} \mathsf{v}_1^{[01]} \mathsf{v}_{2}^{[01]}=\{ab\}\{ba\}\{aa\}, \\
		\mathsf{v}_0^{[10]} \mathsf{v}_1^{[10]} \mathsf{v}_{2}^{[10]}=\{ab\}\{ba\}\{aa\}, \\
		\mathsf{v}_0^{[11]} \mathsf{v}_1^{[11]} \mathsf{v}_{2}^{[11]}=\{ba\}\{ab\}\{ba\}.
		\end{eqnarray*}
		Step 3 and Step 4 then follow from our choice of paths.
	\end{example}
	
	\section{Mixing Properties}
	
	Suppose $X$ is a multidimensional shift space. It is known that shift spaces are associated with various types of dynamical behavior depending on their intrinsic mixing properties; for instance, a shift space $X$ over $\mathbb{Z}^2$ being block-gluing implies that the set of periodic points is dense \cite{BPS-TAMS2010}. In addition to the analogous idea of the aforementioned mixing properties, we consider a collection of mixing properties defined in the sense of the complete prefix code (CPC), from which interesting phenomena also arise in the study of tree-shifts {(cf.~\cite{Ban2017a, Ban2017})}. In this section, the mixing properties for $\mathcal{T}_{}$, $\mathcal{T}_{X}$, and $\mathcal{T}_{A}$, and their connection with $X$ and $\mathcal{T}_{X}$ are presented.
	
	\subsection{Definitions and notations}
	For each word $x \in \Sigma^*$, we say that $x$ is a prefix of $w \in \Sigma^*$, denoted by $x \preceq w$, if there exists $y \in \Sigma^*$ such that $w = xy$, the concatenation of $x$ and $y$. Suppose $x, y \in \Sigma^*$, define the word distance between $x$ and $y$ as
	$$
	d(x, y) := |x| + |y| - 2 \max \{|w|: w \preceq x, w \preceq y\}.
	$$
	From the illustration of Cayley graph of the rooted tree $\Sigma^*$, the word distance $d(x, y)$ reflects the length of the shortest path from $x$ to $y$. Suppose $H, H' \subset \Sigma^*$ are two finite subsets. Define
	$$
	d(H, H') := \min \{d(w, w'): w \in H, w' \in H'\}.
	$$
	Mixing conditions studied in multidimensional shift spaces extend to tree-shifts as follows.
	
	\begin{definition} \label{def:traditional_mixing}
		\item[1.] A tree-shift $\mathcal{T}_{}$ is called \emph{irreducible} (IR) if for any admissible $n$-block $u$ and admissible $m$-block $v$, there exist $w \in \Sigma^\ast \setminus {s(u)}$ and $t \in \mathcal{T}_{}$ such that $t|_{\Delta_n}=u$ and $(\sigma_w t)|_{\Delta_m}=v$.
		
		\item[2.] A tree-shift $\mathcal{T}_{}$ is called \emph{{topologically mixing}\footnote{Note that an analogous definition of  topologically mixing involving the complete prefix code is not mentioned in Definition \ref{def:CPC_mixing}, since it is not quite related to other results discussed in this article.}} (TM) if for any two finite subset of $\Sigma^\ast$, say $H_1$ and $H_2$, there is a positive number $N$ such that for any labeled trees $t_1,t_2 \in \mathcal{T}_{}$ and $w\in \Sigma^\ast$ with $d(w,H_1)\geq N$ there is a labeled tree $t \in \mathcal{T}_{}$ such that  $t|_{H_1}=t_1|_{H_1}$ and $(\sigma_w t)|_{H_2}=t_2|_{H_2}$.
		
		\item[3.] A tree-shift $\mathcal{T}_{}$ is called \emph{block gluing} (BG) if there is a positive number $N$ such that for any admissible $n$-block $u$, admissible $m$-block $v$ and $w\in \Sigma_k^\ast$ with $|w|\geq N+n$ there is a labeled tree $t \in \mathcal{T}_{}$ such that $t|_{\Delta_n}=u$ and $(\sigma_w t)|_{\Delta_m}=v$. 
		
		\item[4.] A tree-shift $\mathcal{T}_{}$ is called \emph{strongly irreducible} (SI) if there is a positive number $N$ such that for any two admissible patterns $u$, $v$ and $w\in \Sigma_k^\ast$ with $d(w,s(u))\geq N$ there is a labeled tree $t \in \mathcal{T}_{}$ such that $t|_{s(u)}=u$ and $(\sigma_w t)|_{s(v)}=v$.
	\end{definition}
	\noindent Note that in the definition above, the patterns are required, as in Section 2.1, to have a finite and prefix-closed support, which yields slightly different definitions of mixing properties compared to the traditional ones on shift spaces over $\mathbb{Z}^d$.
	
	A subset $P$ of $\Sigma^\ast$ is called a \emph{prefix set} if no element in $P$ is a prefix of one another. A prefix set $S \subseteq \Sigma^*$ is called a \emph{complete prefix code} (CPC) if for each $w \in \Sigma_k^\ast$ with $\norm{w} \ge \max_{z \in P} \norm{z}$, there exists $x \in P$ such that $x \preceq w$.
	
	\begin{definition} \label{def:CPC_mixing}
		\item[1.] A tree-shift $\mathcal{T}_{}$ is called \emph{CPC-irreducible} if for any admissible $n$-block $u$ and $m$-block $v$ there is a CPC $P \subseteq \cup_{i \ge n+1} \Sigma^i$ and $t \in \mathcal{T}_{}$ such that $t |_{s(u)}=u$ and $(\sigma_w t)|_{s(v)}=v$ for all $w \in P$.
		\item[2.] A tree-shift $\mathcal{T}_{}$ is called \emph{CPC-block gluing} if there is a CPC $P$ such that for any admissible $n$-block $u$ and admissible $m$-block $v$ there is a labeled tree $t \in \mathcal{T}_{}$ such that $t$ links $u$ and $v$ through $P$.  That is, $t |_{\Delta_n}=u$ and $(\sigma_{w z} t)|_{\Delta_m}=v$ for all $w\in \Sigma^n, z\in P$.
		\item[3.] A tree-shift $\mathcal{T}_{}$ is called \emph{CPC-strongly irreducible} if there is a CPC $P$ such that for any two admissible patterns $u$ and $v$ there is a labeled tree $t \in \mathcal{T}_{}$ such that $t$ links $u$ and $v$ through $P$; i.e., $t |_{s(u)}=u$ and $(\sigma_{w z} t)|_{s(v)}=v$ for all $w\in \partial s(u)$ and $z\in P$, where $\partial s(u) = \{w \in s(u): w \Sigma \cap s(u) = \emptyset\}$ is the boundary of $u$.
		\item[4.] A tree-shift $\mathcal{T}_{}$ is called \emph{uniform CPC-strongly irreducible/block gluing/irreducible} if $\mathcal{T}_{}$ is CPC-strongly irreducible/block gluing/irreducible with $P=\Sigma^n$ for some $n\in \mathbb{N}$.
	\end{definition}
	\subsection{Relations between mixing conditions for tree-shifts $\mathcal{T}_{}$ and $\mathcal{T}_{X}$}
	
	This subsection reveals the relations of mixing conditions between traditional and complete-prefix-code senses. Such relations indicate that tree-shifts are capable of rich dynamical phenomena.
	
	\begin{theorem}  \label{thm:Z_top}
		Suppose $\mathcal{T}_{}$ is a tree-shift. Then,
		\item[1.] $\mathcal{T}_{}$ is SI if $\mathcal{T}_{}$ is CPC USI. 
		\item[2.] $\mathcal{T}_{}$ is BG if $\mathcal{T}_{}$ is CPC UBG. 
		\item[3.] $\mathcal{T}_{}$ is TM if $\mathcal{T}_{}$ is CPC BG. 
		\item[4.] $\mathcal{T}_{}$ is TM if $\mathcal{T}_{}$ is BG. 
	\end{theorem}
	
	\begin{proof}
		We shall prove \textbf{1.}~and \textbf{2.}, since \textbf{3.}~is proved by Ban and Chang in \cite[Proposition~3.6]{Ban2017a} and \textbf{4.}~follow from definition immediately.
		
		\textbf{1.}  Suppose there exists $P=\Sigma^N$ such that every pair of admissible patterns $u, v$ can be linked by a labeled tree $t \in \mathcal{T}_{X}$ through $P$. We show that for every $w \in \Sigma^\ast$ with $d(w, s(u)) \ge N+1$, there exists a corresponding labeled tree $t$ such that $t|_{s(u)} = u$ and $(\sigma_{w} t)|_{s(v)}=v$. Suppose $w=w' w''$, where $w' \in s(u)$ satisfies $\norm{w''} = d(w, s(u)) \ge N+1$. {To utilize CPC USI property, we need to prevent $w$ from having no prefix lying in the boundary $\partial s(u)$ of $u$}. To that end, extend $u$ to an admissible pattern $\overline{u}$ with support $s(\overline{u}) = s(u) \Delta_{\norm{w''}-N} = \{z' z'' \in \Sigma_k^\ast: z' \in s(u), z'' \in \Delta_{\norm{w''}-N}\}$ so that its boundary $\partial s(\overline{u})$ is a CPC and $w = \overline{w}' \overline{w}''$ for some $\overline{w}' \in \partial s(\overline{u})$ and $\norm{\overline{w}''}=N$ {(Note that this extension is always possible since by definition there exists $t \in \mathcal{T}$ such that $t|_{s(u)}=u$).} By CPC USI, there is an labeled tree $t$ with $t|_{s(\overline{u})} = \overline{u}$ and $(\sigma_{z' z''} t)|_{s(v)}=v$ for every $z' \in \partial s(\overline{u})$ and $z'' \in \Sigma^N$. In particular, $(\sigma_w t)|_{s(v)}=(\sigma_{\overline{w}' \overline{w}''} t)|_{s(v)} = v$ since $\overline{w}' \in \partial s(\overline{u})$ and $\overline{w}'' \in \Sigma^N$. 
		
		\textbf{2.} Suppose there exists $P=\Sigma^N$ such that every pair of admissible blocks $u, v$ can be linked by a labeled tree $t \in \mathcal{T}_{}$ through $P$. We claim that for every $w \in \Sigma_k^\ast$ with $d(w, s(u)) \ge N$, there exists a corresponding labeled tree $t$ such that $t|_{s(u)} = u$ and $(\sigma_{w} t)|_{s(v)}=v$. Extend $u$ to an admissible block $\overline{u}$ with support $s(\overline{u}) = \Delta_{\norm{w}-N}$ and $w = \overline{w}' \overline{w}''$ for some $\overline{w}' \in \partial s(\overline{u})$ and $\overline{w}'' \in \Sigma^N$. By CPC UBG, there is an admissible labeled tree $t$ with $t|_{s(u')} = u'$ and $(\sigma_{z' z''} t)|_{s(v)}=v$ for every $z' \in s(\overline{u})$ and $z'' \in \Sigma^N$. In particular, $(\sigma_w t)|_{s(v)}=(\sigma_{\overline{w}' \overline{w}''} t)|_{s(v)} = v$ since $\overline{w}' \in \partial s(\overline{u})$ and $\overline{w}'' \in \Sigma^N$. 
	\end{proof}
	
	{When focusing on the subclass of tree-shifts that comprised of all hom tree-shifts $\mathcal{T}_X$, we obtain Theorem \ref{thm:Z(X)_top} as a refined version of Theorem \ref{thm:Z_top}.}
	
	\begin{theorem} \label{thm:Z(X)_top}
		Suppose $X$ is a shift space. The following implications hold: 
		\item[1.] $\mathcal{T}_{X}$ is CPC USI if and only if $\mathcal{T}_{X}$ is CPC UBG.
		\item[2.] $\mathcal{T}_{X}$ is CPC UBG if and only if $\mathcal{T}_{X}$ is CPC BG.
		\item[3.] $\mathcal{T}_{X}$ is CPC UBG if and only if $\mathcal{T}_{X}$ is BG.
		\item[4.] $\mathcal{T}_{X}$ is CPC IR if and only if $\mathcal{T}_{X}$ is IR.
	\end{theorem} 
	
	\begin{proof}
		\textbf{1.} It suffices to prove the sufficiency. Since $\mathcal{T}_{X}$ is CPC UBG, there exists $N \in \mathbb{N}$ such that $\Sigma^N$ links arbitrary two admissible blocks. We claim that there exists a CPC $P = \Sigma^N$ so that for all admissible patterns $u,v$ there is a labeled tree $t \in \mathcal{T}_{X}$ such that $t|_{s(u)}=u$, $(\sigma_{w z} t)|_{s(v)}=v $ for every $w \in \partial s(u)$ and every $z \in P$. Without loss of generality, $u$ can be extended to some pattern $\overline{u}$ so that $\partial s(\overline{u})$ is a CPC and that $\partial s(u) \subseteq \partial s(\overline{u})$ {by including in $s(\overline{u})$ every element $g \in \Sigma^\ast$ that has all its prefixes lying in $s(u)$}. Suppose $w = w_0 w_1 \cdots w_{n} \in \partial s(\overline{u})$, we write $w^{(i)} = w_0 w_1 \cdots w_i$ for $0 \le i \le n$ and define for every branch $\{w^{(i)}\}_{0 \le i \le n}$ an admissible $n$-block $\overline{u}^{(w)} \in \mathcal{B}_n (\mathcal{T}_{X})$ as follows:
		\begin{equation} \label{eq:uniform_block}
		\overline{u}_g^{(w)} = \overline{u}_{w^{(\norm{g})}}, \forall \  0 \le \norm{g} \le n.
		\end{equation}
		Let $\overline{v}$ be the smallest admissible block such that $v$ is a subpattern of $\overline{v}$. Then, since $\mathcal{T}_{X}$ is CPC UBG, for each $w \in \partial s(\overline{u})$ there is a labeled tree $t^{w} \in \mathcal{T}_{X}$ such that $t^{w}|_{s(\overline{u}^{(w)})}=\overline{u}^{(w)}$ and that $(\sigma_{w z} t^{w})|_{s(\overline{v})}=\overline{v}$ for every $w \in \Sigma^{\norm{w}}$ and every $z \in P$. Define the labeled tree $t$ as
		\begin{equation*}
		t_z:=\begin{cases}
		t^{w}_z, & \text{if } w \preceq z, w \in \partial s(\overline{u}); \\
		u_z, & \text{otherwise. }
		\end{cases}
		\end{equation*}
		It follows that $t_{s(u)}=u$ and $(\sigma_{w z} t)|_{s(v)}=v$ for all $w \in \partial s(u) \subseteq \partial s(\overline{u})$ and all $z \in P$.
		
		\textbf{2.} It is left to prove the sufficiency. Suppose there exists a CPC $P$ such that for every admissible $n$-block $u$ and every admissible block $v$ there is a labeled tree $t$ in which $u$ and $v$ linked through $P$. Let $N=\max \{\norm{u}: u \in P\}$. Without loss of generality, we may assume $0^N \in P$. We claim that for all $u, v \in \mathcal{B} (\mathcal{T}_{X})$ there is a labeled tree $t \in \mathcal{T}_{X}$ such that $t|_{s(u)}=u$ and $(\sigma_{w z} t)|_{s(v)}=v$ for all $w \in \partial s(u)$ and $z \in \Sigma^N$. Since $\mathcal{T}_{X}$ is CPC BG, there is a labeled tree $\hat{t} \in \mathcal{T}_{X}$ such that $\hat{t}|_{s(u)}=u$ and that $(\sigma_{w z} \hat{t})|_{s(v)}=v$ for all $w \in \Sigma^n$ and all $z \in P$. {Let $F:\Sigma^\ast \to \Sigma^\ast$ be the replacement function of branches defined as
		\[
		F_{n,N}(z):=\begin{cases}
		z, & \text{if } \norm{z} \le n; \\
		z_1 z_2 \ldots z_n 0^{\norm{z}-n}, & \text{if } n < \norm{z} \le n+N, \\
		z_1 z_2 \ldots z_n 0^N z_{n+N+1} z_{n+N+2} \ldots z_{\norm{z}}, & \text{if } n+N < \norm{z};
		\end{cases}
		\]
		and define the labeled tree $t$ as $t_z:=\hat{t}_{F_{n,N}(z)}$. Since for every chain $\mathbf{s}$ containing $w \in \Sigma^n$ we have $\pi_{\mathbf{s}}(t)=\pi_{F_{n,N}(\mathbf{s})}(\hat{t})$, it follows that $t \in \mathcal{T}_{X}$ is well-defined and $t_{s(u)}=u$, $(\sigma_{w z} t)|_{s(v)}=v$ for all $w \in \partial s(u)$ and all $z \in \Sigma^N$.} The proof is thus finished.
		
		\textbf{3.} The necessity is part of Theorem \ref{thm:Z_top}. We shall prove the sufficiency. Suppose $N \in \mathbb{N}$ is the minimal gap required to link arbitrary two admissible blocks. We claim that for any admissible $n$-blocks $u$ and admissible block $v$ there is a labeled tree $t$ satisfying $t|_{s(u)} = u$ and $(\sigma_{w z} t)|_{s(v)}=v$ for all $w \in \partial s(u)$ and $z \in \Sigma^N$. For each $w \in \partial s(u)$, define an $n$-block $u^{(w)} \in \mathcal{B}_n (\mathcal{T}_{X})$ as \eqref{eq:uniform_block}.
		Since $\mathcal{T}_{X}$ is BG, for each $w \in \partial s(u)$ there is a labeled tree $t^{w}$ such that $t^{w}|_{s(u^{(w)})}=u^{(w)}$ and that $(\sigma_{w 0^N} t^{w})|_{s(v)}=v$. {Let $F$ be given as above and define the labeled tree $t$ as
		\begin{equation*}
		t_z:=\begin{cases}
		t^{w}_{F_{n,N}(z)}, & \text{if } w \preceq z, w \in \partial s(u); \\
		u_z, & \text{otherwise. }
		\end{cases}
		\end{equation*}
		It follows that $t_{s(u)}=u$ and $(\sigma_{w z} t)|_{s(v)}=v$ for all $w \in \partial s(u)$ and all $z \in P$.
		}
		
		{\textbf{4.} To prove $\mathcal{T}_X$ is CPC IR, we first show that for any $n$-block $u$ with $u_w = u_z$ whenever $\norm{w}=\norm{z}$ and any block $v$, there exist an positive integer $N$ and a labeled tree $t$ such that $t|_{s(u)}=u$ and that $(\sigma_z t)|_{s(v)}=v$ whenever $z \in \Sigma^{n+N}$. If the claim holds, given any arbitrary $n$-block $u$ and block $v$ we can define for every $w \in \partial s(u) = \Sigma^n$ the block $u^{(w)}$ as \eqref{eq:uniform_block}. We then apply our claim to derive an integer $N_w$ and a labeled tree $t^w$ such that $t^w|_{s(u^{(w)})}=u^{(w)}$ and that $(\sigma_z t^w)|_{s(v)}=v$ whenever $z \in \Sigma^{n+N}$. Hence, $u$ and $v$ are linked by a labeled tree $t$ defined as $t_z:=\hat{t}_{F_{n,N_w}(z)}$. More precisely, we have $t|_{s(u)}=u$ and $(\sigma_z t)|_{s(v)}=v$ for every $z \in \cup_{w \in \Sigma^N} w \Sigma^{N_w}$, which is clearly a CPC. We now turn to our claim. Since $\mathcal{T}_X$ is IR, for the $n$-block $u$ with $u_w = u_z$ whenever $\norm{w}=\norm{z}$ and the block $v$ there exist $g \in \Sigma^\ast \setminus s(u) = \cup_{i \ge n+1} \Sigma^i$ and a labeled tree $\hat{t}$ such that $\hat{t}|_{s(u)}=u$ and $(\sigma_g \hat{t})|_{s(v)}=v$. The desired $t$ is then defined as
		\begin{equation*}
		    t_z:=\begin{cases}
		        \hat{t}_{g^{(\norm{z})}} & \text{if } \norm{z} \le \norm{g}; \\
		        \hat{t}_{g z''} & \text{if } z = z' z'', \norm{z'}=\norm{g}, 
		    \end{cases}
		\end{equation*}
		and $N=\norm{g}-n$.}
	\end{proof}
	
	\subsection{Remarks on Theorem \ref{thm:Z_top}}
	
\begin{figure}
		\centering
		\includegraphics[trim=0 0.3cm 0 1.5cm, clip]{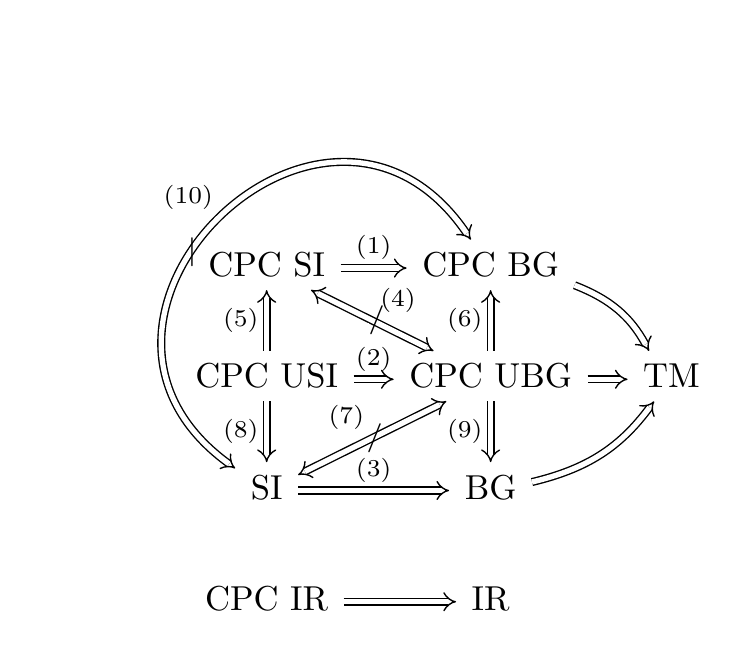}
		\caption{Relations between mixing properties on tree-shifts $\mathcal{T}_{}$}
		\label{fig:mixing_relation_Z}
	\end{figure}

	Figure \ref{fig:mixing_relation_Z} further addresses a completion of Theorem \ref{thm:Z_top}. In this subsection, we give examples or counterexamples for each numbered arrow in the diagram. It remains unknown that whether CPC SI is sufficient for SI/BG, and whether CPC BG is sufficient for BG.
	
	Note that \textbf{(1)}, \textbf{(2)} and \textbf{(3)} are only sufficient but not necessary, which could be shown if we consider the tree-shift $\mathcal{T}_{}$ induced by the forbidden set {$\mathcal{F} = \cup_{n \ge 0} \mathcal{F}_n$ in which $\mathcal{F} _n = \{u \in \{0,1\}^{\Delta_n}: u_w \ne u_z \text{ for some } w,z \in \Sigma^n\}$. Note that every admissible block $u$ or labeled tree $t$ in such a tree-shift can be characterized by its leftmost branch since the symbols within each layer of $\Sigma^\ast$ are unique. In particular, $\mathcal{T}=\{t \in \{0,1\}^{\Sigma^\ast}: t_g = t_w \text{ if } \norm{g}=\norm{w}\}$
	and thus $\pi_{\mathbf{0}}(\mathcal{T})=\{0,1\}^{\mathbb{Z}_+}$.
	}
	
	\textbf{(1)} {To see $\mathcal{T}_{}$ is CPC BG by $\Sigma$, given any two blocks $n$-block $u$ and $m$-block $v$, there is a labeled tree
	\[
	t_g =\begin{cases}
	    u_{0^{\norm{g}}}, & \text{if } 0 \le \norm{g} \le n; \\
	    v_{0^{\norm{g}-n-1}}, & \text{if } n+1 \le \norm{g} \le n+m+1; \\
	    0, & \text{otherwise,} 
	\end{cases}
	\]
	which clearly lies in $\mathcal{T}$ and links $u$ and $v$ through $\Sigma$.} We now show $\mathcal{T}_{}$ is not CPC SI by contradiction. Assume it is CPC SI by a CPC $P$ with $N=\max_{z \in P} \norm{z}$. Let $\overline{w} \in P$ with $\norm{\overline{w}} = N$. Consider admissible patterns $u$, $v$ such that $0 \in \partial s(u) \cap \Sigma$ and that there exists $\overline{z} \in s(u) \cap \Sigma^{N+1}$ with $u_{\overline{z}} \ne v_\epsilon$. If $\mathcal{T}_{}$ is CPC SI, then there is a labeled tree $t$ such that $t|_{s(u)} = u$ and $(\sigma_{w z} t)|_{s(v)} = v$ for every $w \in \partial s(u)$ and $z \in P$. In particular, $t_{\overline{z}} = u_{\overline{z}} \ne v_\epsilon = t_{0 \overline{z}}$, while $\norm{\overline{z}} = \norm{0 \overline{w}} = N+1$, contradicting the definition of $\mathcal{F}$. 
	
	\textbf{(2)} {We have shown that $\mathcal{T}_{}$ is CPC UBG in \textbf{(1)}, and it cannot be CPC USI by definition, since it is demonstrated above that $\mathcal{T}_{}$ is not CPC SI.}
	
	\textbf{(3)} We prove that $\mathcal{T}_{}$ is not SI by contradiction, assuming that $\mathcal{T}_{}$ is SI by an integer $N$. Let $\overline{w} \in \Sigma^{N}$ be fixed. Consider admissible patterns $u$, $v$ such that $s(u)=\{\epsilon, 0\} \cup \{1^i\}_{N+1 \ge i \ge 1}$ and $u_{1^{N+1}} \ne v_\epsilon$. By our assumption of SI, there is a labeled tree $t$ such that $t|_{s(u)} = u$ and $(\sigma_{0 \overline{w}} t)|_{s(v)} = v$. In particular, $t_{1^{N+1}} = u_{1^{N+1}} \ne v_\epsilon = t_{0 \overline{w}}$ while $\norm{1^{N+1}} = \norm{0 \overline{w}} = N+1$, contradicting the hypothesis of $\mathcal{F}$. 
	
	On the other hand, that a tree-shift is CPC SI is neither sufficient nor necessary for it to be CPC UBG ({see} \textbf{(4)}). The reader may refer to \cite[Example~3.12]{Ban2017a} for a detailed proof. 
	
	Furthermore, that a tree-shift is SI is neither sufficient nor necessary for the tree-shift to be CPC BG ({see} \textbf{(10)}). If we combine this with the results above, it yields that SI is neither necessary nor sufficient for CPC UBG ({see} \textbf{(7)}), and \textbf{(8)} and \textbf{(9)} are only sufficient but not necessary. Non-necessity part is already proved in the above case. As for the non-sufficiency, consider $\mathcal{T}_{} \subseteq \{0,1\}^{\Sigma^\ast}$ to be induced by the forbidden set $\mathcal{F} = \{(1;1,1),(0;1,1)\}$. Then, it is SI by a distance 2, since $\Sigma^2$ always leaves sufficient vacancy for balance in symbol 0 and symbol 1. However, one can show that it is not CPC BG by contradiction, assuming there exists a CPC $P$ through which all admissible patterns $u, v$ can be linked. Since $P$ is a CPC, there exists $z \in \Sigma^\ast$ such that $z0, z1 \in P$. Now we consider $u=1$ and $v=1$, for which a labeled tree $t \in \mathcal{T}_{}$ satisfies $t_\epsilon=1$ and $t_{w}=1$ for all $w \in P$. In particular, $(t_z; t_{z 0}, t_{z 1}) \in \mathcal{F}$, which contradicts the assumption.
	
	\subsection{Relations between $X$ and $\mathcal{T}_{X}$}
	
	Aside from investigating relations between mixing conditions for $\mathcal{T}$ and $\mathcal{T}_X$, respectively, it is of interest to study the relations of mixing conditions between a shift space $X$ and its corresponding hom tree-shift $\mathcal{T}_X$. This subsection reveals that the implications are only one-directional.
	
	\begin{theorem} \label{thm:Z(X)_top_with_X}
		Suppose $X$ is a shift space. The following implications hold: 
		\item[1.] $X$ is mixing if $\mathcal{T}_{X}$ is TM.
		\item[2.] $X$ is transitive if $\mathcal{T}_{X}$ is IR.
	\end{theorem}
	
	\begin{proof}
	    \textbf{1.} Suppose $u = u_0 \cdots u_n$ and $v = v_0 \cdots v_n$ are admissible blocks in $X$. Let $\overline{u} \in \mathcal{B}_n (\mathcal{T}_{X})$ and $\overline{v} \in \mathcal{B}_m (\mathcal{T}_{X})$ be constructed as $\overline{u}_g=u_{\norm{g}}$ for $0 \le \norm{g} \le n$ and $\overline{v}_g=v_{\norm{g}}$ for $\norm{g} \le m$. Since $\mathcal{T}_{X}$ is TM, there is an $N=N(s(u), s(v)) \in \mathbb{N}$ such that for all $w \in \partial s(\overline{u})$ and $z \in \Sigma^{N+i}$, $i \ge 0$, there is a $t \in \mathcal{T}_{X}$ satisfying $t|_{s(\overline{u})}=\overline{u}$ and $(\sigma_{w z} t)|_{s(\overline{v})}=\overline{v}$. By choosing a chain $\mathbf{s}$ containing $w z$ and defining $x=\pi_{\mathbf{s}}(t) \in X$, it is clear $x_{[0,n]}=u$ and $x_{[n+\norm{z},n+m+\norm{z}]}=v$. This implies $X$ is mixing.
	    
	    {\textbf{2.} Suppose $u = u_0 \cdots u_n$ and $v = v_0 \cdots v_n$ with $\overline{u}$ and $\overline{v}$ defined as above. Since $\mathcal{T}_{X}$ is IR, there is a $w \notin s(\overline{u}) = \Delta_n$ and a $t \in \mathcal{T}_{X}$ satisfying $t|_{s(\overline{u})}=\overline{u}$ and that $(\sigma_{w} t)|_{s(\overline{v})}=\overline{v}$. By choosing a chain $\mathbf{s}$ containing $w$ and defining $x=\pi_{\mathbf{s}}(t) \in X$, it is clear that $x_{[0,n]}=u$ and $x_{[\norm{w},\norm{w}+m]}=v$, where $\norm{w} > n$. This implies $X$ is transitive.}
	\end{proof}
	
	{Before we present the example illustrating the converse of Theorem \ref{thm:Z(X)_top_with_X} could fail, a related class of shift spaces called \emph{bounded density shifts} (see \cite{stanley2013bounded}) is introduced. Given a function $f: \mathbb{Z}_+ \to [0,\infty)$, the bounded density shift associated with $f$ is defined as
	\begin{equation*}
	    \Psi_f:=\left\{x \in \mathbb{Z}_+^{\mathbb{Z}}: \sum_{r=i}^{i+p-1} x_r \le f(p), \forall p \in \mathbb{Z}_+,i \in \mathbb{Z}\right\},
	\end{equation*}
	whose one-sided version is given as
	\begin{equation*}
	    \Psi_f^+ = \left\{x_{[0,\infty)}: x \in \Psi_f\right\}.
	\end{equation*}
	It is clear from the definition that a bounded density shift $\Psi_f$ (and $\Psi_f^+$ as well) can be associated with infinitely many $f$; in other words, there are always infinitely many distinct functions $g: \mathbb{Z}_+ \to [0,\infty)$ satisfies $\Psi_f=\Psi_g$ (or $\Psi_f^+ = \Psi_g^+$ for the one-sided case). Nevertheless, as discussed in \cite[Lemma 2.1]{stanley2013bounded} there exists a \emph{canonical} $f$ such that if $\Psi_f=\Psi_g$ (or $\Psi_f^+=\Psi_g^+$ for the one-sided case), then $f \le g$. This requirement of $f$ makes it an integer-valued function.}
	
	{
	Bounded density shifts play an essential role in the following discussion.
	}
	
	\begin{remark} \label{rmk:Z(X)_top}
		Suppose $X$ is a shift space and $\mathcal{T}_X$ is its corresponding hom tree-shift. Figure \ref{fig:mixing_relation_Z(X)} illustrates Theorems \ref{thm:Z(X)_top} and \ref{thm:Z(X)_top_with_X}.
		\begin{figure}
			\centering
			\includegraphics{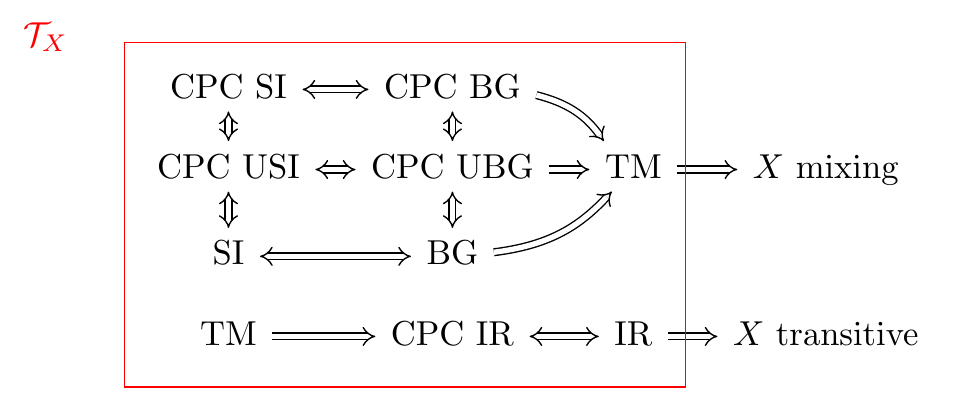}
			\caption{Relations between mixing properties on tree-shifts $\mathcal{T}_{X}$ and shift spaces $X$}
			\label{fig:mixing_relation_Z(X)}
		\end{figure}
		Notably, the converse of Theorem \ref{thm:Z(X)_top_with_X} may not hold in general. For instance, one may consider that even shift $X$ has mixing property while the hom tree-shift $\mathcal{T}_{X}$ is not TM. Suppose $X$ is the even shift, which is known to be mixing. However, $\mathcal{T}_{X}$ is not even CPC IR by considering admissible blocks
		\begin{equation} \label{eq:counterexample_even_shift}
		u=1, \hspace{2em} v=\vcenter{\hbox{\includegraphics[]{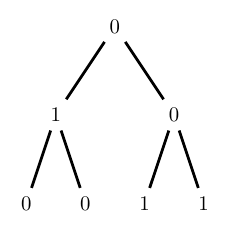}}},
		\end{equation}
		{for if it is CPC IR, there exist $t \in \mathcal{T}_{X}$ and $g \in \Sigma^\ast$ such that $t_\epsilon=1$ and $(\sigma_g t)|_{s(v)}=v$. In particular, $t|_{g0}=t|_{g10}=1$. By choosing two chains $\mathbf{s}^1$, $\mathbf{s}^2$ such that $g0 \in \mathbf{s}^1$ and $g10 \in \mathbf{s}^2$, we observe that $\pi_{\mathbf{s}^1}(t)_{[0,\norm{g}+2]} = \pi_{\mathbf{s}^1}(t)_{[0,\norm{g}]} 10$ and $\pi_{\mathbf{s}^2}(t)_{[0,\norm{g}+2]} = \pi_{\mathbf{s}^1}(t)_{[0,\norm{g}]} 01$. Since $\pi_{\mathbf{s}^1}(t)_{[0,\norm{g}]} \ne 0^{\norm{g}+1}$, either $\pi_{\mathbf{s}^1}(t)_{[0,\norm{g}+2]}$ or $\pi_{\mathbf{s}^2}(t)_{[0,\norm{g}+2]}$ is a forbidden block of the even shift $X$ and we reach a contradiction that $t \notin \mathcal{T}_X$.}
	    
		Another remark can be made here is that CPC UBG is not necessary for TM, for which an example is given as follows. {Let $X=\Psi_f^+$ be a bounded density shift with $f$ defined as 
		\[
		f(n) := \lceil \log_3 (n+1) \rceil
		\]
		for $n \ge 0$. It is shown in \cite[Lemma 2.13]{stanley2013bounded} that $f$ is canonical for $\Psi_f$ and unboundedness of $f$ implies $\Psi_f$ is mixing by \cite[Theorem 2.14]{stanley2013bounded}. Consequently, $\Psi_f^+$ is also mixing.} We then show that $\mathcal{T}_{X}$ is TM by proving given two arbitrary shapes $H_1, H_2$ and $M=\max_{z \in H_1 \cup H_2} \norm{z}$, {there exists an integer
		\begin{equation} \label{eq:TM_criterion}
		    N = \max_{1 \le \ell \le M+1} \{ \min \{N_\ell \in \mathbb{N}: f(\ell+N_\ell) \ge f(M+1) + f(\ell)\} \}
		\end{equation}
		as mentioned in Definition \ref{def:traditional_mixing}.} More specifically, if $H_1, H_2$ and $N$ are given as above, then for all $\overline{w} \in \Sigma^\ast$ satisfies $d(\overline{w}, H_1) \ge N$ and all $t_1, t_2 \in \mathcal{T}_{}$, one can define the following labeled tree:
		\[
		t_z := \begin{cases}
		(t_1)_z, & \text{ if } z \in H_1; \\
		(t_2)_{\overline{w}'}, & \text{ if } z = \overline{w} \overline{w}', \overline{w}', \in H_2; \\
		0, & \text{ otherwise.}
		\end{cases}
		\]
		It suffices to show that $t \in \mathcal{T}_{X}$, i.e. $\pi_{\mathbf{s}} (t) \in X$ for every chain $\mathbf{s}$. Suppose $\mathbf{s}$ is an arbitrary chain that contains $\overline{w}$. Without loss of generality, the set $H_1 \cap \mathbf{s}$ and $w H_2 \cap \mathbf{s}$ can be expressed as follows:
		\begin{align*}
		H_1 \cap \mathbf{s} = \{\overline{z} s_1 s_2 \cdots s_i\}_{0 \le i \le n},\\
		\overline{w} H_2 \cap \mathbf{s} = \{\overline{w} s_1' s_2' \cdots s_i'\}_{0 \le i \le m},
		\end{align*}
		where $\overline{z} \le \overline{w}$ and $\norm{\overline{w}}-\norm{\overline{z}} \ge N$. It is sufficient to prove that $\pi_{\mathbf{s}} (t) \in \mathcal{T}_{X}$. {To be more precise, for every $\norm{\overline{z}} \le i \le \norm{\overline{z}}+n$ and $\norm{\overline{w}} \le j \le \norm{\overline{w}}+m$, equation \eqref{eq:TM_criterion} together with $n,m \le M+1$ and monotonicity of $f$ implies that
		\begin{equation*}
		    \begin{aligned}
		        \sum_{r=i}^{j}  \pi_{\mathbf{s}}(t)_r &= \sum_{r=i}^{\norm{\overline{z}}+n}  \pi_{\mathbf{s}}(t)_r + \sum_{r=\norm{\overline{w}}}^{j}  \pi_{\mathbf{s}}(t)_r \\
		        &\le f(\norm{\overline{z}}+n-i+1) + f(M+1) \\
		        &\le f(\norm{\overline{z}}+n-i+1+N) \\
		        &\le f(p-i).
		    \end{aligned}
		\end{equation*}
		This shows that $\pi_\mathbf{s}(t) \in X$ for every chain $\mathbf{s}$ and thus $t \in \mathcal{T}_X$. 
		} 
		
		{It is left to show $\mathcal{T}_{X}$ is not CPC UBG. Otherwise, suppose $\mathcal{T}_{X}$ is CPC UBG by $P=\Sigma^N$. Note that \cite[Corollary 2.6 and Lemma 2.7]{stanley2013bounded} assures that 
		\[
		x=(x_i)_{i \in \mathbb{Z}_+}=\left(f(i+1) - f(i)\right)_{i \in \mathbb{Z}_+}
		\]
		lies in $X$. Since the alphabet of $X$ consists of only symbol 0 and 1, it follows from the definition of $f$ that if $y \in X$ satisfies $y_{[0,3^N-1]} = x_{[0,3^N-1]}$ then $y_{i}=0$ for every $3^N \le i \le 3^{N+1}-2$. To see this, note that for $3^N \le n \le 3^{N+1}-1$,
		\begin{equation*}
		        N+1 = f(n) \ge \sum_{i=0}^{n-1} y_i \ge \sum_{i=0}^{3^N-1} y_i = f(3^N) = N+1.
		\end{equation*}
		Hence, if we consider the admissible block $(3^{N}-1)$-block $u$ defined by $u_g=x_{\norm{g}}$ and $0$-block $v=1$ and suppose $u$ and $v$ are linked by some $t \in \mathcal{T}_X$, then $t_g = 1$ for every $\norm{g} = 3^{N}+N-1$ by the definition of CPC UBG. However, $t_g = \pi_{\mathbf{0}}(t)_{3^{N}+N-1}=1$ with $3^N \le 3^{N}+N-1 \le 3^{N+1}-1$ implies that $\pi_{\mathbf{0}}(t) \notin X$, which is a contradiction.}
		
		Since $\mathcal{T}_{X}$ is not UBG if $X$ is the bounded density shift as given above, it is neither CPC BG nor BG by Theorem \ref{thm:Z(X)_top}.
	\end{remark}
	
	{
	\begin{remark} \label{rmk:axial_product}
        It is noteworthy that the example regarding \eqref{eq:counterexample_even_shift} also distinguishes the hom tree-shifts from the axial product on semigroups. In fact, if $X$ is the even shift and $\mathcal{T}_X'$ is its axial power on $\Sigma^\ast$, then $u$ and $v$ can be linked by an $t \in \mathcal{T}_X'$, as opposed to $\mathcal{T}_X$. More explicitly, define $t$ as
		\begin{equation*}
		    t_g:=\begin{cases}
		        1, & \text{if } g \in \{\epsilon, 000,0010,0011\}; \\
		        0, & \text{otherwise.}
		    \end{cases}
		\end{equation*}
	    Then, $t \in \mathcal{T}_X'$, $t|_{s(u)}=u$, and $(\sigma_0 t)|_{s(v)}=v$.
	\end{remark}
	}
	
	Whenever $X = X_A$ is a Markov shift determined by the binary matrix $A$, the relations of mixing conditions between $X_A$ and $\mathcal{T}_A$ are equivalent as characterized by Theorem \ref{thm:Z(X_M)_top}.
	
	\begin{theorem} \label{thm:Z(X_M)_top}
		Suppose $A$ is an adjacency matrix. The following implications hold:
		\item[1.] $X_A$ is mixing if and only if $\mathcal{T}_{A}$ is CPC UBG.
		\item[2.] $X_A$ is transitive if and only if $\mathcal{T}_{A}$ is IR.
	\end{theorem} 
	\begin{proof}
		For simplicity, for the rest of the proof, we denote by $[w]$ the set of all prefixes of $w \in \Sigma^\ast$ with $\epsilon$ included. 
		
		\textbf{1.} The sufficiency is part of Theorem \ref{thm:Z(X)_top}. It remains to show the necessity. Since $X_A$ is a Markov shift, there is an integer $N \in \mathbb{N}$ such that for all admissible blocks $\hat{u},\hat{v} \in \mathcal{B}(X_A)$ and all $n \ge N$, there is a configuration $x \in X_A$ such that $x_{[0, \norm{\hat{u}}-1]}=\hat{u}$ and $x_{[\norm{\hat{u}}+n-1,\norm{\hat{u}}+n+\norm{\hat{v}}-2]}=\hat{v}$. We use this property to show that for every admissible $m$-blocks $u$ and every admissible block $v$ in $\mathcal{T}_{A}$, there is a labeled tree $t$ such that $t|_{s(u)}=u$ and $(\sigma_{w z} t)|_{s(v)}=v$ for all $w \in \Sigma^m$ and $z \in \Sigma^N$. Suppose $\overline{z} \in \partial s(v)$ is fixed. By mixing property of $X_A$, for every branch $[\overline{w}]$ with $\overline{w} \in \partial s(u)$ there exists $x^{(\overline{w})} \in X_A$ such that $x^{(\overline{w})}|_{[0, \norm{\overline{w}}-1]}=u_{[\overline{w}]}$ and that $x^{(\overline{w})}|_{[\norm{\overline{w}}+N-1, \norm{\overline{w}}+N+\norm{\overline{z}}-2]} = v_{[\overline{z}]}$. Take, by definition of admissible block, labeled trees $t^u, t^v \in \mathcal{T}_{A}$ such that $t^u|_{s(u)}=u$ and $t^v|_{s(v)}=v$ and define the desired labeled tree
		\begin{equation*}
		t_z:=\begin{cases}
		t^v_{z'}, & \text{if } z = \overline{w} \overline{z} z', \overline{w} \in \partial s(u); \\
		x^{(\overline{w})}_{\norm{z'}}, & \text{if } z = \overline{w} z', \overline{w} \in \partial s(u), \norm{z'} \le N; \\
		t^u_{z}, & \text{otherwise}, \\
		\end{cases}
		\end{equation*}
		which lies in $\mathcal{T}_{X}$ and satisfies $t|_{s(u)}=u$ and $(\sigma_{w z} t)|_{s(v)}=v$.
		
		{\textbf{2.} Since the sufficiency is proved in Theorem \ref{thm:Z(X)_top}, we prove the necessity. Let $u$, $v$, $t^u$, $t^v$ have the same meaning as above. Since $X_A$ is irreducible, for every branch $[\overline{w}]$ with $\overline{w} \in \partial s(u)$ there exists $x^{(\overline{w})} \in X_A$ such that $x^{(\overline{w})}|_{[0, \norm{\overline{w}}-1]}=u_{[\overline{w}]}$ and that $x^{(\overline{w})}|_{[\norm{\overline{w}}+N_{\overline{w}}-1, \norm{\overline{w}}+N_{\overline{w}}+\norm{\overline{z}}-2]} = v_{[\overline{z}]}$ for some $N_{\overline{w}}$ depending on $\overline{w}$. In a similar manner, the desired labeled tree linking $u$ and $v$ is defined as
		\begin{equation*}
		t_z:=\begin{cases}
		t^v_{z'}, & \text{if } z = \overline{w} \overline{z} z', \overline{w} \in \partial s(u); \\
		x^{(\overline{w})}_{\norm{z'}}, & \text{if } z = \overline{w} z', \overline{w} \in \partial s(u), \norm{z'} \le N_{\overline{w}}; \\
		t^u_{z}, & \text{otherwise}, \\
		\end{cases}
		\end{equation*}
		and the proof is completed.}
	\end{proof}
	
	\begin{remark} \label{rmk:Z(M)_top}
		Suppose $A$ is an adjacency matrix. We see from Theorem \ref{thm:Z(X)_top} that CPC UBG implies TM, from Theorem \ref{thm:Z(X)_top_with_X} that TM implies $X_A$ is mixing, and from Theorem \ref{thm:Z(X_M)_top} that $X_A$ being mixing implies $\mathcal{T}_{A}$ is CPC UBG. Hence, all of the above are equivalent. See Figure \ref{fig:mixing_relation_Z(M)}.
		\begin{figure*}
			\centering
			\includegraphics[]{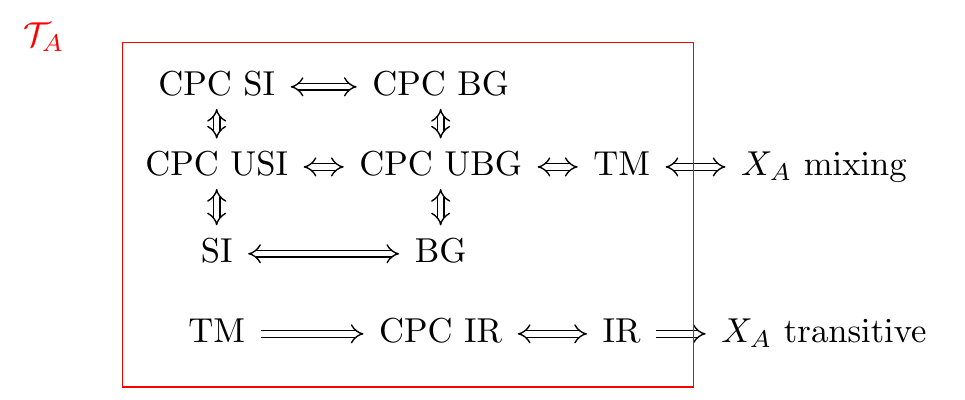}
			\caption{Relations between mixing properties on hom tree-SFTs $\mathcal{T}_{A}$ and Markov shifts $X_A$}
			\label{fig:mixing_relation_Z(M)}
		\end{figure*}
	\end{remark}
	\section*{Acknowledgement}
	We would like to sincerely thank the anonymous referee for providing inspiring comments and helpful suggestions for the first draft of this article. These significantly improve the readability and solidify the validity of the theorems in the paper.
	\bibliographystyle{amsplain}
	\bibliography{reference}
\end{document}